\documentclass[a4paper, twoside, 11pt]{smfart}
\usepackage[latin1]{inputenc}
\evensidemargin \oddsidemargin 

\usepackage{amsmath,amssymb,amsthm}
\usepackage[francais,english]{babel} 
\usepackage[arrow, matrix]{xy}
\usepackage{joe2002,joethm}
\usepackage{polyzetas}

\def\vdeux#1#2{#2}
\numberwithin{equation}{section}

\title
{Doubles mélanges des polylogarithmes multiples aux racines de l'unité}

\author{Georges Racinet}
\address{Mathematisches Institut, Einsteinstra\ss e 62, D--48149, Münster} 
\email{racinet@math.uni-muenster.de}
\urladdr{http://www.dma.ens.fr/~racinet}
\date{\today}

\begin{document}
\frontmatter
\selectlanguage{french}
\begin{abstract}
   Les valeurs des fonctions polyzêtas aux entiers strictement positifs 
fournissent une solution au système 
d'équations des associateurs de Drinfel'd, aux nombreuses applications en
algèbre quantique. Vues comme intégrales itérées, ce sont les périodes 
du groupoïde fondamental motivique de $\mathbb{P}^1\setminus\{0,1,\infty\}$, 
d'où un système fondamental de relations algébriques, 
qui implique celui des associateurs mais n'est pas explicite.  

  On étudie ici la combinatoire d'un autre système de relations, 
les \emph{doubles mélanges}, qui provient de manipulations élémentaires de 
séries et d'intégrales. On montre qu'il partage une propriété importante 
avec les associateurs et les relations \og{}motiviques\fg{}, est conséquence de
ces dernières et définit une algèbre de polynômes sur $\mathbb{Q}$ 
(théorème d'Écalle).  On obtient ces résultats pour les nombres plus généraux
que sont les polylogarithmes multiples aux racines de l'unité de Goncharov.
\end{abstract} 
\begin{altabstract}
  The values at positive integers of the polyzeta functions are solutions 
of the polynomial equations arising from Drinfeld's associators, which have
numerous applications in quantum algebra. Considered as iterated integrals 
they become periods of the motivic fundamental groupoid of 
$\mathbb{P}^1\setminus\{0,1,\infty\}$. From there comes a fundamental, yet 
no more explicit, system of algebraic relations; it implies the system of 
associators. 

   We focus here on the combinatorics properties of another system of 
relations, the ``double shuffles'', which comes from elementary series and
integrals manipulations. We show that it shares an important property with
associators and ``motivic'' relations, is implied by the latter and defines
a polynomial algebra over $\mathbb{Q}$ (Écalle's theorem). We obtain these
results for more general numbers: values of Goncharov's multiple polylogarithms
at roots of unity.
\end{altabstract}	
\maketitle
\section{Introduction}\label{sec_intro}
 
Les valeurs aux entiers strictement positifs des fonctions polyzêtas 
\begin{equation}
  \z(s_1,\ldots,s_r)\ \ass\ 
\sum_{n_1>n_2>\cdots>n_r>0} 1/{n_1^{s_1}n_2^{s_2}\cdots n_r^{s_r}}
\end{equation}
fournissent des solutions réelles de problèmes de nature purement algébrique,
souvent dans le domaine de l'algèbre quantique. 
Déjà considérées par Euler dans le cas $r=2$, elles sont portées à l'attention
du public par Zagier \cite{ZagierECM} au moment où elles  
apparaissent comme les 
coefficients de l'associateur de Drinfel'd $\Phi_\KZ$ \cite{DrinQTQH, LeMur}.
Plus récemment, Kontsevitch les a retrouvées dans son isomorphisme de 
formalité \cite{Kon2}.   

Dans une série d'articles  \cite{Gonch98, GonchDuke, Gonch2001}, Goncharov
s'intéresse à des objets plus généraux, les polylogarithmes multiples. Ce sont
les fonctions de $r$ variables
	complexes définies dans le polydisque unité par:
\begin{equation}\label{def_L}
  L_{s_1,\ldots,s_r}(z_1,\ldots,z_r)\ \ass\ 
\sum_{n_1>n_2>\cdots>n_r>0} \frac{z_1^{n_1}z_2^{n_2}\cdots z_r^{n_r}}
{n_1^{s_1}n_2^{s_2}\cdots n_r^{s_r}},
\end{equation} 
où $s_1,\ldots,s_r$ sont des entiers strictement positifs. Ces séries 
entières sont de rayon de convergence $1$ et ne divergent sur le 
polycercle unité que pour $(s_1,z_1)=(1,1)$. Les polyzêtas ne sont que les
valeurs de ces fonctions lorsque les $z_i$ valent tous 1. 

Les valeurs des polylogarithmes multiples sont aussi celles d'intégrales 
itérées du type
\begin{equation}\label{def_I}
I_{[0,1]}(a_1, \ldots, a_p) =
\int_{0\leq t_p\leq\cdots\leq t_1\leq 1} \bigwedge_{i=1}^p \omega_{a_i}(t_i),
\end{equation}
avec $\omega_a(t)=adt/(a-t)$ pour $a\neq 0$ et $\omega_0(t)=dt/t$. 

Dans le cas des polyzêtas, c'est à Kontsevitch que l'on doit la
formule reliant ces deux types d'objets. 
En convenant que $0^k$ désigne la séquence 
formée du chiffre $0$ répété $k$ fois, elle s'écrit:
\begin{equation}
\z(s_1,\ldots,s_r) = I_{[0,1]}(0^{s_1-1},1,0^{s_2-1},1,\ldots,0^{s_r-1},1),
\end{equation} 
et Goncharov l'a généralisée aux polylogarithmes multiples.

\subsection{Relations de mélange, exemples}
Notre propos est la combinatoire des relations de \og{}double mélange\fg{} 
des valeurs des 
polylogarithmes multiples. Pour des raisons de finitude, on se limitera
au cas où les $z_i$ 
\vdeux{sont des racines $n$\emes{} l'unité.}{
parcourent le groupe $\mub_n(\CM)$ des racines $n$\emes{} l'unité, $n$ étant 
fixé.} 
Ces relations 
proviennent de manipulations formelles élémentaires, sur les écritures
(\ref{def_L}) et (\ref{def_I}). 

Par les séries entière, on obtient par exemple
\begin{multline}\label{exemple_rel2}
L_{s_1}(z_1)L_{s_2}(z_2) = \sum_{n_1,n_2>0} \frac{z_1^{n_1}z_2^{n_2}}
{n_1^{s_1}n_2^{s_2}} = \left(\sum_{n_1>n_2>0} + \sum_{n_2>n_1>0} + 
\sum_{n_1=n_2>0}\right) \frac{z_1^{n_1}z_2^{n_2}}{n_1^{s_1}n_2^{s_2}}\\
 = L_{s_1,s_2}(z_1,z_2)
+ L_{s_2,s_1}(z_2,z_1) + L_{s_1+s_2}(z_1z_2),
\end{multline}
\vdeux{}{et une formule plus compliquée pour le produit de deux séries du 
type (\ref{def_L}).}

D'un autre coté, 
le produit de deux simplexes se décompose en union de simplexes, ce qui fournit
la relation de mélange des intégrales itérées:
\begin{equation}\label{exemple_relsh}
 I_{[0,1]}(a_1,\ldots,a_p)I_{[0,1]}(b_1,\ldots b_q)
= \sum_{\s\in\SG_{p,q}} I_{[0,1]}(a_{\s^{-1}(1)}, \ldots, a_{\s^{-1}(i)})
\end{equation}
Dans cette formule, $\SG_{p,q}$ désigne l'ensemble des 
$(p,q)$-battages\footnote{En anglais, les {\em shuffles} ; \og{}mélange\fg{}
est un terme moins précis.}: les
permutations de $\{1,\ldots,p+q\}$ qui sont croissantes sur $\{1,\ldots,p\}$ 
et $\{p+1,\ldots, q\}$.

On a donc deux systèmes de relations dont la combinatoire semble similaire.
On leur adjoindra la relation de régularisation, qui provient 
de l'annulation formelle de certaines divergences. 
Les arguments 
qu'on utilise pour l'obtenir sont empruntés à Boutet de Monvel \cite{Boutet}.
On essaie d'en préciser l'historique. 

On forme ainsi le système DMR (Doubles Mélanges \& Régularisation).   
Une conjecture de Kontsevitch et Zagier prévoit 
\vdeux{que le système DMR }{qu'il est \emph{complet}, \ie{} qu'il }engendre toutes les relations algébriques entre polyzêtas.
Les valeurs aux racines de l'unité satisfont également à d'autres relations,
distribution et poids un, qu'on regroupe dans le système DMRD, dont 
on sait qu'il n'est pas complet en général. 

Au cours de la section \ref{sec_descr}, on donne une description duale 
du système DMRD:
 
Étant donné un sous-groupe fini multiplicatif $\G$ de $\CM^*$, on considère
un alphabet $\XB_\G$ dont les éléments, notés $x_\s$, sont indexés par $\GZ$.
les relations DMRD s'expriment de manière compacte sur la série 
\vdeux{}{génératrice}{} non commutative
\begin{equation}
   \ItC = \sum_{p\in\NM, a_1,\ldots,a_p\in \GZ}
I_\g(a_1,\ldots,a_p)x_{a_1}\cdots x_{a_p},
\end{equation}
qui n'est autre, lorsque $\G=\{1\}$, que le $\Phi_\KZ$ de Drinfel'd.%
\footnote{Nos notations et celles de \cite[\S 5]{DrinQTQH} se correspondent
par $x_0=A, x_1=-B$.} 
La définition de $\ItC$ nécessite de \og{}régulariser\fg{} les intégrales 
itérées divergentes, c'est-à-dire de leur donner un sens convenable. 
Cela sera fait directement sur les séries\vdeux{.}{ génératrices.} 

On a choisi d'établir directement les relations DMRD dans cette optique duale,
cela n'étant guère plus long qu'un rappel de la présentation habituelle et 
la preuve de l'équivalence. Cela nous a naturellement conduits
à détailler les démonstrations de nombre de propriétés élémentaires. 

\subsection{Aspects motiviques}
Les intégrales itérées du type (\ref{def_I}) ont un sens motivique qui 
explique une grande partie des propriétés algébriques des polyzêtas. On en
donne ici un bref aper\c cu, renvoyant à \cite{DelPi1, Gonch2001, Del01}
pour un exposé complet.

Soit $X=\Punrac{n}$, vu comme schéma sur un corps cyclotomique 
$F$ de degré $n$ et choisissons un plongement complexe $\s$ de $F$. 
Le groupoïde fondamental motivique de $X$, défini par 
Goncharov \cite{Gonch2001}, est un
pro-objet de la catégorie tannakienne des motifs de Tate mixtes 
sur $\OC_{F,S}$ de {\em loc. cit}. 
Il a les images décrites par Deligne dans \cite{DelPi1} 
par les foncteurs fibres de réalisation $\omega_{B,\s}$
(de Betti, relativement à $\s$) et $\omega_\dr$ (de de Rham).
Notamment, en réalisation de de Rham (définie sur $F$)
$\pi_1^\dr(X;a,b)$ est indépendant des points-base $a$ et $b$ et
s'identifie aux exponentielles de Lie 
\vdeux{de $\XB$}{formées sur $\XB\ \ass\ \XB_{\mu_n(\CM)}$}.
\vdeux{D'après \cite{Del01}, on }{On }a 
une $\QM$-structure canonique $\omega_\gr$ sur 
\vdeux{la réalisation de de Rham. }
{le foncteur-fibre $\omega_\dr$ \cite{Del01}}.

La série $\ItC$ s'interprète comme 
l'image du chemin $[0,1]$ par l'isomorphisme de comparaison
\begin{equation}\label{eq_BDR}
 \xymatrix@1{\pi_1^{B,\s}(X;(0,1), (1,0))(\CM)\ar[r]&
\pi_1^\dr(X; (0,1), (1,0))(\CM)}, 
\end{equation}
où $(0,1)$ et $(1,0)$ sont les points-base tangentiels de Deligne. 
Par la théorie générale des catégories tannakiennes, les isomorphismes 
entre les foncteurs-fibre $\omega_B$ et $\omega_\gr$ forment une variété 
pro-algébrique affine sur $\QM$, torseur sous le groupe $\GC_\gr$ des 
automorphismes de $\omega_\gr$. Ce groupe admet une décomposition 
$\GC_\gr = \Gm \ltimes \UC_\gr$, où $\UC_\gr$ est le noyau de l'action 
sur le motif de Tate $\QM(1)$. Par définition, $\GC_\gr$ et $\UC_\gr$ 
commutent à toute flèche de la forme $\omega_\gr(f)$. 

Soient $G$ et $H$ deux éléments de $\sernc{\k}{\XB}$. 
On étudie en \ref{subsec_MT} le produit $G\pmt H$, défini si le terme constant
de $G$ est 1. Ce produit fait des séries à terme constant 1 de 
$\sernc{\k}{\XB}$ un groupe pro-unipotent, noté $\MT(\k)$, 
qui agit sur $\sernc{\k}{\XB}$. 
\vdeux{}{Cette opération apparaît, sous diverses formes, chez Drinfel'd, 
Goncharov et Ihara.} 

\def\mot{\text{\sf M}}
Il est expliqué en détail dans \cite{Del01} comment l'action de 
\vdeux{$\UC_\dr$}{$\UC_\gr$} sur 
$\pi_1^\gr(X; (0,1), (1,0))$ se factorise par $\MT$, autrement dit, se fait
par la loi $\pmt$.  Deligne utilise 
l'image de $\UC_\gr$ dans $\MT$ pour construire une sous-variété fermée 
$\mot$ de $\AM^1\times\MT$ dont 
 $\ItC$ est un point complexe,
exprimant ainsi les relations entre valeurs aux racines de l'unité 
des polylogarithmes multiples qui sont conséquences 
de la comparaison Betti-de Rham.\footnote{%
Cette construction nous est plus facile d'accès que celle, essentiellement 
équivalente,
des polyzêtas motiviques $\z_\MC$ de Goncharov \cite{GonchICM}.} 
Une variante de la conjecture des périodes
de Grothendieck prévoit que ce système est complet.
On a une flèche naturelle $\mot\to\AM^1$. Le tout satisfait à la propriété
(M) ci-dessous:

\begin{defi}
   Une sous-variété fermée $V$ de $\pi_1^\dr(X; (0,1), (1,0))$ munie 
   d'une flèche $V\to\AM^1$ dont on note $V_\l$ les fibres au-dessus 
   de $\l\in\AM^1(\k)=\k$ a la propriété (M) si $V_0$ est un sous-schéma en 
   groupes de $\MT$ et $V\to\AM^1$ est un torseur trivial sous 
   \vdeux{$\MT$.}{l'action de $V_0$ par $\pmt$.}  
\end{defi}
On ne connait pas de description de $\mot$ par un système explicite 
d'équations. Dans le cas des polyzêtas, 
\vdeux{}{\ie{} $X=\Puntrois$,}{}
la variété $\Ass$ des associateurs 
de Drinfel'd  \cite{DrinQTQH}\footnote{Drinfel'd la note $\MC$. Elle est 
en fait construite dans $\AM^1\times\pi_1^\dr(X; (0,1), (1,0))$ et
munie de la première projection. Ce n'est pas une différence essentielle.} 
est un candidat. Les équations la définissant 
sont liées aux conditions de cohérence de Mac-Lane 
et Drinfel'd prouve que $\Ass$ a la propriété (M) 
\cite[props. 5.5 \& 5.9]{DrinQTQH}. 
La nature motivique de $\Ass$ fait peu de doute, apparaissant
clairement dans la présentation de Bar-Natan \cite{BarNatGT}.%
\footnote{On peut interpréter
ses catégories $\widehat{PAB}$ et $\widehat{PACD}$ comme 
les \og{}tours\fg{} Betti et de Rham des groupoïdes fondamentaux des 
espaces de configuration de $n$ points
sur la droite affine, restreints à un jeu de points-base tangentiels.} 
Le groupe 
$\Ass_0$ est la variante de de Rham du groupe de Grothendieck-Teichmüller ;
Drinfel'd le note $\GRT_1$. 
L'étude de
$\Ass$, et particulièrement la recherche d'associateurs rationnels est 
en soi un problème important, du fait des nombreuses applications
\cite{Tam2, KEbial1, BarNatNAT, PCinv}. 

Le versant $\ell$-adique de ces objets a été étudié indépendamment par Ihara 
\cite{Ihcont,IhFest,IhIsrael} ; Hain et Matsumoto \cite{HainMat} ont récemment 
accompli des progrès dans cette direction. 

\subsection{Résultat principal}
Dans la section \ref{sec_formel}, on définit, par leurs $\k$-points, des
éléments de $\sernc{\k}{\XB}$,   
les variétés pro-algébriques $\DMRP$ et $\DMRD$ correspondant aux 
systèmes DMR et DMRD. D'après la section \ref{sec_descr}, $\ItC$ est 
un point complexe de $\DMRP$ et $\DMRD$.  

Une première propriété de l'espace tangent $\dmr$ à $\DMRP$ au voisinage 
de la solution évidente $1$ amène à considérer une flèche $\DMRP\to\AM^1$,
attribuant à une série un de ses coefficients. Pour la série $\ItC$, 
ce coefficient n'est autre que $\z(2)$, si l'on travaille avec au plus 
deux racines de l'unité ; dans les autres cas, c'est un multiple rationnel de 
$2i\pi$.  

On définit alors $\MT$ et on en étudie les propriétés les plus immédiates,
notamment différentielles. 
Cela permet enuite d'expliciter l'énonce de notre résultat principal\vdeux{}{,
annoncé dans \cite{crasjoe1,crasjoe2}}: 
$\DMRP$ et $\DMRD$ 
\vdeux{}{ont la propriété (M).} 

\subsection{Preuve}
La démonstration se fait en deux étapes. D'abord, on prouve que $\dmr_0$ 
et $\dmrd_0$ sont des sous-algèbres de Lie de 
$\mt$\vdeux{}{, l'algèbre de Lie de $\MT$, }
et que les $\DMRP_\l$ 
sont stables par l'action de l'exponentielle, au sens 
de $\MT$, de $\dmr_0$. 
C'est l'objet de la section \ref{sec_tangente}. La démonstration 
repose sur l'enchaînement de plusieurs miracles combinatoires, et l'auteur
doit bien avouer ne pas avoir d'argument plus direct à proposer: la nature
motivique des relations du type (\ref{exemple_rel2}) n'est pas claire.

Au cours de la section \ref{sec_transitivite}, on complète la démonstration du 
théorème en montrant que l'action de $\exp(\dmr_0)(\k)$ sur $\DMRP_\l(\k)$ est 
transitive pour tout $\k$. Cela résulte de la nature d'espace tangent de 
$\dmr_0$ et de l'existence d'une solution particulière: $\ItC$. L'argument 
est inspiré de \cite{BarNatGT}. Au prix d'une certaine gymnastique formelle,
essentiellement tautologique, on le développe dans un cadre un peu plus 
général. 

On étudie ensuite 
quelques conséquences directes du théorème. Notamment, on exhibe  
certains éléments irréductibles de $\dmrd_0$ égaux, dans le cas des 
polyzêtas, à ceux de $\grt_1$ définis par Drinfel'd \cite{DrinQTQH}. 
Dans tous les cas, ils engendrent l'algèbre
de Lie de $\mot_0$ dans la description de Deligne \cite{Del01}.
Autrement dit, les relations motiviques impliquent les
relations DMRD. 

\subsection{Remerciements}Les plus vifs sont adressés à Pierre Cartier,
qui a encadré ce travail, et à Pierre Deligne pour ses patientes 
explications. La présentation de la section \ref{sec_descr} profite 
directement des raccourcis de \cite{DelNotes}. 

Merci également à S. Aicardi, J. Bellaïche, L. Boutet de Monvel, F. Digne,
B. Enriquez, 
A. Goncharov, I. Marin, J.--C. Novelli, M. Petitot et C. Reutenauer.

Ce travail a été effectué alors que l'auteur préparait une thèse à
l'université de Picardie-Jules-Verne et était membre du département 
de mathématiques de l'École Normale Supérieure. L'auteur 
bénéficie actuellement d'une bourse post-doctorale de l'ESF à 
l'université de Münster, où la rédaction a été achevée.  

\subsection{Terminologie, notations, rappels, abus}\label{subsec_conv}
Les anneaux sont supposés commutatifs et unifères, les algèbres sont
associatives et unifères; les cogèbres sont coassociatives et coünifères.
Un $\QM$-anneau est un anneau contenant $\QM$. Le groupe à un 
élément est noté $\triv$. On fixe une clôture algébrique $\CM$ de $\RM$.

On adopte le point de vue fonctoriel sur les schémas, utilisé notamment 
par Demazure et Gabriel \cite{Demaz}: 
un schéma affine 
sur $\QM$ est un foncteur représentable de la catégorie des $\QM$-anneaux.
Une variété pro-algébrique affine est une limite projective de schémas 
algébriques affines. C'est encore un schéma affine. 
La lettre $\k$ désigne en général un $\QM$-anneau quelconque, qui sert 
de variable pour ces foncteurs. 

On utilise implicitement le produit tensoriel gradué et le produit 
tensoriel complété. Ce dernier est parfois noté $\whot$ si l'on tient 
à préciser. 

Un $\QM$-espace vectoriel gradué $V$ dont les composantes homogènes
sont de dimension finie donne lieu à un schéma vectoriel 
$\k\mapsto V\whot\k$. On le notera simplement $V$, utilisant $V(\k)$ si 
nécessaire. 
La $\QM$-algèbre correspondante est l'algèbre symétrique du dual gradué 
de $V$ (voir \cite[p. 147]{Demaz} pour le cas usuel). 

On note $\mglk(C)$ l'ensemble des éléments {\em\glks{}}\footnote{Cela
a été jugé préférable à l'utilisation répétée de \og{}group-like\fg{}.}
 d'une cogèbre 
$(C,\eps,\D)$, c'est-à-dire les $x\in C$ qui vérifient 
$\D x = x\ot x$ et $\eps(x)=1$. 
Dans le cas où $C$ est une bigèbre graduée, rappelons que l'application 
exponentielle est une bijection
de l'ensemble des éléments primitifs du complété $\wh{C}$ dans l'ensemble
des éléments \glks{} de $\wh{C}$. 

Pour un ensemble $Z$, on désigne respectivement 
par $\assl{\k}{Z}$ et $\liel{\k}{Z}$
l'algèbre associative libre et l'algèbre de Lie libre 
sur $Z$ à coefficients dans $\k$. Dans tous les cas que l'on considérera, 
une graduation 
convenable permettra de définir les algèbres filtrées complètes
correspondantes, notées $\sernc{\k}{Z}$ et $\whliel{\k}{Z}$ et de traiter
leurs éléments comme des séries formelles. \vdeux{}{On renvoie au livre 
de Reutenauer \cite{Reut} pour un exposé et une bibliographie complets 
à propos des algèbres de Lie libres.}   

Pour ces objets libres, on qualifiera de {\em substitution} un morphisme dans
la catégorie appropriée défini par son action sur les générateurs. 

Pour $w$ mot en $Z$ et $\Phi\in\sernc{\k}{Z}$, la notation 
$\scal{\Phi}{w}$ désigne le coefficient de $w$ dans $\Phi$ et est étendue
à tout élément de $\assl{\k}{Z}$ par linéarité sur le membre de droite. 
On l'utilisera également pour les éléments de puissances tensorielles.
On en fait l'usage le plus limité possible.

Une $\QM$-algèbre de Lie $\gG$ est considérée comme incluse dans son 
algèbre enveloppante universelle $\Ue{}{\gG}$. 
C'est une bigèbre de Hopf pour le coproduit $\D$ pour lequel les 
éléments de $\gG$ sont primitifs, la coünité $\eps$ par 
$\gG\subset\ker\eps$ et l'antipode $S$ par $S(x)=-x$ pour tout $x\in\gG$.  
Si $\gG$ est graduée et ses composantes homogènes de dimension finie, 
on note $\Uec{\k}{\gG}$ le complété à coefficients dans $\k$, munie de 
la structure de bigèbre obtenue par l'extension continue des scalaires 
$\whot\k$. 

On utilisera fréquemment la propriété immédiate suivante: 
un morphisme d'algèbres (resp. une 
dérivation) entre deux 
bigèbres est un morphisme de cogèbres (resp. une codérivation) 
si l'identité appropriée est vraie sur un système de générateurs. 
Par exemple, si $\gG$ est stable pour une dérivation de $\U{\gG}$, celle-ci
est une codérivation.

\section{Description combinatoire des relations}\label{sec_descr}   
 \subsection{Divergences logarithmiques}\label{subsec_divlog}
Notre but est ici de décrire les 
étapes de Boutet de Monvel qui nous permettront 
d'établir la relation de régularisation. 

On note $D$ le disque ouvert unité de $\CM$ et $\HC(D)$ l'algèbre 
des fonctions holomorphes sur $D$. La lettre $t$ désigne une variable 
formelle.

\begin{defi}
  Soit $\divlogs$ l'ensemble des suites $(S_N)_{N>0}$ de nombres 
complexes admettant  un développement asymptotique du type
$$ S_N = \As((S_N))(\log(N)) + o(\log^\alpha(N)/N), $$ 
avec $\As((S_N))\in\CM[t]$ et $\alpha\in\NM$. 

De même, on notera $\divlogf$ l'ensemble des fonctions  $f$ de $\HC(D)$ 
admettant, lorsque $z$ tend vers 1 dans $D$, un développement du type
$$ f(z) = \As(f)(\log(1-z)) + o\left(\aug(-\log(1-z))^\alpha\right) $$
avec $\As(f)\in\CM[t]$ et $\alpha\in\NM$. 
\end{defi}
Il est clair que $\divlogs$ et $\divlogf$ sont des sous-$\CM$-algèbres de 
$\CM^\NM$ et $\HC(D)$, respectivement ; que les polynômes $\As((S_N))$ et 
$\As(f)$ sont uniquement déterminés par ces conditions et que les applications 
$\As$ ainsi définies sont toutes deux des morphismes d'algèbres, grâce
à la nature du reste. L'ambiguïté de notation n'est pas gênante. 
\sss{}\label{sss_serfonc}
À une suite $(S_N)_{N>0}$, faisons correspondre 
l'unique série entière 
$\sum_{n>0}u_n z^n$ telle que $S_N = \sum_{n=1}^N u_n$.  

\begin{prop}\label{prop_kerAs} 
Si $(S_N)_{N>0}$ appartient à $\divlogs$ et $\As((S_N))=0$, la série entière 
associée définit une fonction $f$ de $\divlogf$, vérifiant $\As(f)=0$.
\end{prop}
\begin{proof}
On suppose donc qu'on a un entier $\alpha$  tel que 
$$S_N\ \ssa\ \sum_{n=0}^{N-1} u_n = o(log^\alpha(N)/N)$$
Effectuons une transformation d'Abel:
$$f(z)=\sum_{n>0}u_nz^n= \sum_{N>0} 
\left(\sum_{n=1}^{N}u_n\right)(z^N-z^{N+1}) = 
(1-z)\sum_{N>0} z^N o(\log^{\alpha}(N)/N)$$
Il suffit donc d'obtenir l'existence de $\beta>0$ tel que 
\begin{equation}
\sum_{N>0}z^N\log^{\alpha}(N)/N = o((-\log(1-z)^\beta) 
\end{equation}
Découpons la somme en deux, autour de $1/(1-z)$. En majorant $z$ par 1, on 
obtient
$$ \sum_{0<N<(1-z)^{-1}}z^N\log^\alpha N/N \leq \int_{1}^{(1-z)^{-1}}dt\log^\alpha t/t
= (-\log(1-z))^{\alpha+1}/(\alpha+1)$$
La fonction $\log^\alpha(t)/t$ est décroissante au voisinage
de $+\infty$. Pour $z$ suffisamment proche de 1, et 
$N\geq (1-z)^{-1}$, 
on a donc $\log^\alpha(N)/N\leq(-\log(1-z))^\alpha(1-z)$, dont on déduit 
$$\sum_{N\geq (1-z)^{-1}}z^N\log^{\alpha}(N)/N \leq 
(-\log(1-z))^{\alpha}(1-z)\sum_{N\geq (1-z)^{-1}}z^N \leq  
z^{1/1-z}(-\log(1-z))^{\alpha}
$$
Comme $z^{1/1-z}$ tend vers $\exp(-1)$, ceci donne le
résultat, avec $\beta=\alpha+1$.
\end{proof}
\sss{}\label{sss_comp}
La série harmonique $H_N=\sum_{n=1}^N1/n$ appartient à $\divlogs$.
On a $\As(H_N)=\gamma+t$, où $\gamma$ est la constante d'Euler. Le morphisme 
d'algèbres $\As : \divlogs\to\CM[t]$ est donc surjectif. De même, en 
considérant la fonction $L_1 : z\mapsto -\log(1-z)$, on voit que 
$\As : \divlogf\to\CM[t]$ est également surjectif. 

Soit $k$ un entier positif et 
$(H_{k,N})_{N>0}$ la suite associée à la fonction $(L_1)^k/k!$. 
Il sera prouvé au paragraphe \ref{sss_aposteriori} 
que $(H_{k,N})_{N>0}$ appartient à 
$\divlogs$ et que le polynôme $\As((H_{k,N})_{N>0})$ est de 
degré $k$. Pour tout $(S_N)_{N>0}\in\divlogs$, il existe donc une 
combinaison linéaire de $(S_N)$ et des $(H_{k,N})$ qui est dans le noyau 
de $\As$. On en déduit, grâce à la proposition \ref{prop_kerAs}, le résultat
ci-dessous, qui ne servira pas avant \ref{sss_asympit}: 
\begin{cor}\label{cor_comp}
La fonction associée à un élément de $\divlogs$ appartient à $\divlogf$. 
Il existe une application $\CM$-linéaire $\comp$ faisant commuter le 
diagrame ci-dessous.
$$\xymatrix{\divlogs\ar[d]\ar[r]^-{\As}&\CM[t]\ar[d]^{\comp} \\
	    \divlogf\ar[r]^-{\As}&\CM[t]} $$ 
\end{cor}
\sss{Itération}\label{sss_divlogLN} 
Considérons les sommes partielles associées aux polylogarithmes 
multiples:
\begin{equation}\label{eq_defLN}
L^N_{s_1,\ldots,s_r}(z_1,\ldots,z_r)\ \ass\ 
\sum_{N>n_1>n_2>\cdots>n_r>0} \frac{z_1^{n_1}z_2^{n_2}\cdots z_r^{n_r}}
{n_1^{s_1}n_2^{s_2}\cdots n_r^{s_r}}
\end{equation}
\begin{prop}
   Soit $s$ un entier strictement positif, $z$ dans le disque unité fermé
de $\CM$ et $(S_N)_{N>0} \in \divlogs$, 
La suite $(T_N)_{N>0}$ 
$$ T_N\ \ass\ \sum_{0<n<N} z^n S_n/n^s$$ 
est encore dans $\divlogs$. La dérivée du polynôme $\As(T_N)$ est $\As(S_N)$. 
La suite $T_N$ est convergente si $(s,z)\neq (1,1)$. 
\end{prop}
\begin{preuve}
On a donc un polynôme $P$ tel que $S_n = P(\log(n)) + o(\log^\alpha(n)/n)$
Dans un premier temps, on majore $|z|$ par 1. 

Dans le cas $s\geq 2$, la convergence de $T_N$ est évidente. 
Il faut donc évaluer le reste de la série. 
Au pire, on a à faire à $\sum_{n>N} \log^{\beta}(n)/n^2$, où $\beta$ est
le degré de $P$. Ceci est équivalent à l'intégrale
$$I_\beta\ \ass\ \int_N^{+\infty} \log^{\beta}x^{-2}dx, $$
qui vérifie $I_\beta = \beta\log^{\beta-1}(N)/N+\beta I_{\beta-1}$ avec 
$I_0=1/N$.  
     
Si $s=1$, on a un terme de la forme $o(\log^\alpha(n))/n^2$, qui 
est traité comme ci-dessus, et le terme $P(\log(n))/n)$. 
Il suffit de traiter le cas $P(t)=t^k$, pour tout $k\in\NM$. 

En mettant $\log(n+1)-\log(n)$ en facteur dans 
$\log^{k+1}(n+1)-\log^{k+1}(n)$, on obtient:
\begin{equation}\label{eq_conv}
\log^k(n)/n = (\log^{k+1}(n+1)-\log^{k+1}(n))/(k+1) + o(\log^k(n)/n^2)
\end{equation}
La sommation du premier terme donne $\log^{k+1}(N)/(k+1)$. Le reste
est traité comme précédemment. Ceci prouve les deux premières assertions. 

Comme $S_n/n$ tend vers 0, la convergence est dans le cas $s=1, z\neq1$
une application classique de la transformation d'Abel.
%
\end{preuve}
La suite $L^N_{s_1,\ldots,s_r}(z_1,\ldots,z_r)$ se déduit
de $L^N_{s_2,\ldots,s_r}(z_2,\ldots,z_r)$ comme ci-dessus.
\begin{cor}\label{prop_LNdivlog}
Pour tous entiers strictement positifs $s_1,\ldots,s_r$ 
et tous $z_1,\ldots,z_r$ de module au plus 1, la suite
$(L^N_{s_1,\ldots,s_r}(z_1,\ldots,z_r))_{N>0}$ appartient à $\divlogs$. 
Elle est convergente si et seulement si $(s_1,z_1)\neq(1,1)$.
\end{cor}
\sss{Remarques} Boutet de Monvel considère des développements plus fins,
du type $\sum_{i\geq 0}P_i(\log N)/N^i$. Il aboutit à la jolie formule
$\comp = (d/dt)!$, où la factorielle désigne le 
développement en série entière de la fonction $\Gamma$ 
d'Euler au voisinage de $1$. 

Dans \cite{thesejoe}, on utilisait 
des restes en $o(N^{-\alpha-1})$, avec $\alpha\in\RM^*_+$,
plus grossiers qu'ici. Les démonstrations pour les restes en 
$\log^\alpha(N)/N$ sont celles de \cite{DelNotes}.    
\Subsection{Relations de mélange des intégrales itérées} 
\label{subsec_intit}
\begin{defi} 
Soient $a_1,\ldots,a_p\in\CM$ et $\gamma : [0,1]\to\CM$ un chemin. 
On considère l'intégrale éventuellement divergente:
$$ I_\gamma(a_1,\ldots,a_p)\ \ass\ 
\int_0^1 (\g^*\omega_{a_1})(t_1)
\int_0^{t_1} (\g^*\omega_{a_2})(t_2)
\cdots \int_0^{t_p} (\g^*\omega_{a_p})(t_p),$$
avec, pour tout $a\in\CM:$ 
$$\omega_a(t)\ \ass\ \left\{\begin{array}{ccc}
dt/(a^{-1}-t)&\text{si}&a\neq0\\
dt/t&\text{si}&s=0\end{array}\right.
$$
\end{defi}
Lorsque $\gamma$ évite les singularités $a_1,\ldots,a_p$, cette intégrale
est bien définie et ne dépend que de la classe d'homotopie de $\g$ dans 
$\CM\setminus\{a_1,\ldots,a_p\}$. 
\vdeux{}{Si $p=0$, elle vaudra $1$ par convention.} 

Pour $(a,b)\in\CM^2$, on notera $[a,b]$ 
le chemin  $t\mapsto a + t(b-a)$. 


Goncharov utilise une généralisation de la formule de 
Kontsevitch. En convenant
que $0^k$ désigne la séquence formée du chiffre $0$ répété $k$ fois, elle 
s'écrit:    
\begin{prop}\label{prop_formKont}
 Soient $r\in\NM, s_1,\ldots,s_r\in\NM^*$ et $z_1,\ldots,z_r$ 
des nombres complexes tels que $0<|z_i|\leq 1$. Pour tout $z$ du disque 
unité ouvert de $\CM$, on a
$$ L_{s_1,\ldots,s_r}(z_1,\ldots,z_{r-1},z_{r}z) = I_{[0,z]}
\left(0^{s_1-1}, z_1, 0^{s_2-1}, z_1z_2, \ldots, 0^{s_r-1}, 
z_1z_2\cdots z_r\right)$$ 
\end{prop}
\sss{Séries de Chen}\label{def_Xb}
Soit $\G$ un groupe commutatif fini, noté multiplicativement.
On considère un alphabet $\XB_\G=\{(x_\s)\}$ indexé par $\GZ$ 
Lorsqu'il n'y a pas d'ambiguïté sur $\G$, 
on écrira simplement $\XB$ pour $\XB_\G$. 

On appellera {\em poids} la graduation qui attribue à chaque élément
de $\XB$ le degré 1. Rappelons que le coproduit $\D$ de $\assl{\QM}{\XB}$, 
vue comme bigèbre enveloppante de $\liel{}{\XB}$, est 
homogène et est donné par $ \D x_\s = 1\ot x_\s + x_\s\ot1,\ 
\text{pour tout}\ \s\in\G\cup\{0\}$.  

Supposons $\G$ inclus dans $\CM$. 
À tout chemin $\g$ sur $\PunG$, on associe l'élément suivant de 
$\sernc{\CM}{\XB}$, appelé \emph{série de Chen}. On peut la voir comme une 
sorte de série génératrice non-commutative des intégrales itérées. 
$$ \ItC_\g\ \ass\ \sum_{p\in\NM, a_1,\ldots,a_p\in \G\cup\{0\}}
I_\g(a_1,\ldots,a_p)x_{\s_1}\cdots x_{\s_p}
$$

\begin{prop}[Relation de mélange]\label{prop_Itabdiag}
Pour tout chemin $\g$ sur $\PunG$,  la série $\ItC_\g$ est \glke{} dans 
$(\sernc{\CM}{\XB},\pcdot,\D)$.
\end{prop}
C'est un fait bien connu, qui remonte aux annees 50 \cite{Chen57}, 
et il en existe de nombreuses démonstrations.
Une première méthode \cite{Ree58}, 
est d'établir la formule (\ref{exemple_relsh}), qui 
permet d'interpréter $I_\gamma$ comme un morphisme d'algèbres à valeurs dans 
$\CM$, depuis le dual gradué de $\assl{\QM}{\XB}$ muni du produit de battage 
$\sh$; ce dernier est dual du coproduit $\D$ \cite{Reut}. 

On peut aussi voir $\ItC_{[a,z]}$ comme solution de l'équation différentielle 
\begin{equation}\label{eq_diff}
   d\ItC_{[a,z]} = \Omega\ItC_{[a,z]}\qtext{avec}
\Omega\ \ass\ \sum_{\s\in\GZ} \omega_\s(z)x_\s
\end{equation}
La primitivité de $\Omega$ implique que $\D\ItC_{[a,z]}$ et 
$\ItC_{[a,z]}\ot\ItC_{[a,z]}$ sont solutions de $d\JC = \D(\Omega)\JC$, avec
la même condition initiale $\JC(a)=1\ot 1$. 

Plus généralement, $\ItC_\g$ est le transport parallèle \vdeux{}{de $1$} 
le long de $\g$ de la connexion intégrable 
\vdeux{$d+\Omega$}{$d+(\text{multiplication à gauche 
par $\Omega$})$} sur le fibré trivial
$\PunG\times \sernc{\CM}{\XB}$. C'est ainsi qu'est 
explicitée la comparaison Betti-de Rham pour $\pi_1(\Punrac{n})$
\cite[12.16]{DelPi1}. Pour $\G=\triv$, 
l'équation (\ref{eq_diff})
est la réduction à $\Puntrois$ du système de Knizhnik-Zamolodchikov $\KZ_3$ 
utilisée par Drinfel'd \cite{DrinQTQH}.

\sss{Convergence partielle}\label{sss_convpart}
Pour $(n,\nu)\in\NM^*\times\G$, soit $y_{n,\nu}=x_0^{n-1}x_\nu$. Notons 
$\YB_\G$ l'ensemble des $y_{n,\nu}$, ou simplement $\YB$ s'il n'y a pas 
d'ambiguïté sur $\G$. La sous-algèbre de $\assl{\QM}{\XB}$ 
engendrée par $\YB$ est clairement libre sur $\YB$. Comme espace vectoriel,
elle est engendrée par les mots en $\XB$ ne se finissant pas par $x_0$ et
il est pratique de l'identifier au quotient de $\assl{\QM}{\XB}$ par l'idéal 
à droite homogène $I_0=\assl{\QM}{\XB}x_0$. On notera 
$\pi_Y$ la projection correspondante. 
 
Comme $x_0$ est primitif,  l'idéal $I_0$ est un coïdéal pour $\D$. 
Ceci fait du quotient $\assl{\QM}{\XB}/I_0$ une cogèbre graduée, 
qu'on notera par abus $(\assl{\QM}{\YB},\D)$. 
Par homogénéité, ces constructions passent à l'extension des scalaires 
complétée. Ceci ne fait pas de $(\assl{\QM}{\YB}, \pcdot, \D)$ une 
bigèbre: le coproduit-quotient $\D$ n'est plus un morphisme d'algèbres. 

Considérons $I_{[a,z]}$ comme une fonction de la variable $a$ de 
$]0,1[$. Par la proposition \ref{prop_Itabdiag}, 
L'élément $\pi_Y(\ItC_{[a,z]})$ de $(\sernc{\CM^{]0,1[}}{\YB},\D)$
est \glk{} pour $\D$, car $\pi_Y$ est un morphisme de cogèbres. Ses 
coefficients, du type $I_{[a,z]}(a_1,\ldots,a_r)$ avec $a_r\neq 0$, admettent 
une limite lorsque $a$ tend vers $0$ dans $]0,1[$ 
(prop. \ref{prop_formKont}). Le passage à la limite étant un morphisme 
d'algèbres, on en déduit: 
\begin{prop}\label{prop_Itzdiag}
 Pour tout $z\in]0,1[$, la série 
$$\ItC_{[0,z]} = \sum_{p\in\NM, a_1,\ldots,a_{p-1}\in\GZ, a_{p}\in\G}
I_{[0,z]}(a_1,\ldots,a_p)x_{a_1}\cdots x_{a_p}$$
est un élément \glk{} de la 
cogèbre $(\sernc{\CM}{\YB},\D)$ 
\end{prop}
\sss{}
Soient $\ps$ et $\qs$ (produits et quotients successifs) 
les automorphismes lin{\'e}aires, inverses l'un de l'autre 
de $\assl{\QM}{\XB}$ donnés par 
\begin{eqnarray}
\ps(x_0^{n_1}x_{\s_1}x_0^{n_2}x_{\s_2}\cdots x_0^{n_r}x_{\s_r}x_0^{n_{r+1}}) 
&=& x_0^{n_1}x_{\s_1}x_0^{n_2}x_{\s_1\s_2}\cdots 
x_0^{n_r}x_{\s_1\s_2\s_r}x_0^{n_{r+1}} \\
\qs(x_0^{n_1}x_{\s_1}x_0^{n_2}x_{\s_2}\cdots x_0^{n_r}x_{\s_r}x_0^{n_{r+1}}) 
&=& x_0^{n_1}x_{\s_1}x_0^{n_2}x_{\s_2\s_1^{-1}}\cdots x_0^{n_r}x_{\s_r\s_{r-1}^{-1}}x_0^{n_{r+1}},
\end{eqnarray}
où les $n_i$ et les $\s_i$ sont quelconques dans $\NM$ et dans $\G$. Il est
clair que $\ps$ et $\qs$ sont homogènes pour le poids, ce qui permet d'étendre
leurs définitions à $\sernc{\k}{\XB}$ et que leurs restrictions 
à $\sernc{\k}{\YB}$, qui est stable, sont données par:
\begin{eqnarray}
\ps(y_{s_1,\s_1}y_{s_2,\s_2}\cdots y_{s_r,\s_r}) &=&
y_{s_1,\s_1}y_{s_2,\s_1\s_2} \cdots y_{s_r,\s_1\s_2\cdots\s_r} \\
\qs(y_{s_1,\s_1}y_{s_2,\s_2}\cdots y_{s_r,\s_r}) &=&
y_{s_1,\s_1}y_{s_2,\s_2\s_1^{-1}}\cdots y_{s_r,\s_r\s_{r-1}^{-1}},
\end{eqnarray}
pour tous entiers $s_1,\ldots,s_r$ strictement positifs.
 
Considérons également dans la série génératrice 
\begin{equation}
\LCC_z\ \ass\ \sum_{s_1,\ldots,s_r, \s_1,\ldots,\s_r} 
L_{s_1,\ldots,\s_r}(\s_1,\ldots,z\s_r)y_{s_1,\s_1}\cdots y_{s_1,\s_r}
\end{equation}
Avec ces notations, la proposition \ref{prop_formKont} se réécrit:
\begin{prop} \label{prop_kontserz}
On a $\LCC_z=\qs(\ItC_z)$ dans $\sernc{\CM^{]0,1[}}{\YB}$
\end{prop}
\sss{Convergence totale}
Le sous-espace vectoriel $\assl{\QM}{\XB}_\cv$ de $\assl{\QM}{\XB}$ 
engendré par les mots en $\XB$
ne se finissant pas par $x_0$ et ne commen\c cant pas par $x_1$ est lui-aussi
une sous-algèbre graduée, qu'on identifiera au quotient 
$\assl{\QM}{\XB}/I$ avec $I=(x_1\assl{\QM}{\XB} + \assl{\QM}{\XB}x_0)$. 
On notera $\pi_\cv$ la projection correspondante.
Par primitivité de $x_0$ et $x_1$, le sous-espace $I$ est un coïdéal, 
d'où un coproduit-quotient, toujours noté $\D$ sur $\assl{\QM}{\XB}_\cv$. 

La proposition \ref{prop_kontserz} montre que 
les coefficients de $\pi_\cv(\ItC_z)$ sont les fonctions 
associées comme en \ref{sss_serfonc} aux suites 
$(L^N_{s_1,\ldots,s_r}(\s_1,\ldots,\s_r))_{N>0}$, qui sont convergentes (cor.
\ref{prop_LNdivlog}). Le lemme d'Abel permet de répéter le même argument
qu'en \ref{sss_convpart}:
\begin{prop}\label{prop_Itcvdiag}
La série $\ItC_\cv\ \ass\ \lim\limits_{a\to0^+, b\to1^-} 
\pi_\cv(\ItC_a^b)$ est un élément \glk{} de $(\sernc{\CM}{\XB}_\cv,\D)$. 
\end{prop}

L'espace vectoriel $\assl{\QM}{\XB}_\cv$ est également la sous-algèbre de
$\assl{\QM}{\YB}$ engendrée par les mots en $Y$ ne commen\c cant pas par 
$y_{1,1}=x_1$. On l'identifie encore 
 au quotient $\assl{\QM}{\YB}/y_{1,1}\assl{\QM}{\YB}$. Cela nous amènera à
noter indifféremment $\assl{\QM}{\YB}_\cv$ pour $\assl{\QM}{\XB}_\cv$. De
plus, on voit facilement que $\ps$ et $\qs$ commutent à $\pi_\cv$.  

\begin{prop} \label{prop_kontsercv}
La série génératrice 
\begin{equation}
   \LCC_\cv\ \ass\ \sum_{s_1,\ldots,s_r, \s_1,\ldots,\s_r, 
(s_1,\s_1)\neq(1,1)} 
L_{s_1,\ldots,\s_r}(\s_1,\ldots,\s_r)y_{s_1,\s_1}\cdots y_{s_1,\s_r}
\end{equation}
est égale à $\qs(\ItC_\cv)$ dans $\sernc{\CM}{\YB}_\cv$
\end{prop}
En effet, $\LCC_\cv$ et $\qs(\ItC_\cv)$ sont les limites terme à terme 
de $\pi_\cv(\LCC_z)$ et $\qs(\pi_\cv(\ItC_{[0,z]}))$ lorsque $z$ tend vers $1^-$ 

\subsection{Relation de mélange des sommes itérées}
On se propose ici de décrire de manière compacte toutes les relations du type
(\ref{exemple_rel2}), lorsque les variables sont dans un sous-groupe 
multiplicatif fini $\G$ fixé de $\CM$.   
Comme en \ref{subsec_intit}, cela se traduira par la
\glkte{} de certaines séries génératrices. 
\sss{Définitions} \label{sss_defYDet}
Dans $\assl{\QM}{\YB}$, on étend la notation $y_{n,\nu}$ au cas $n=0$ par la
convention  
\begin{equation}\label{conv_y0}
 y_{0,\s}\ \ass\ \left\{\begin{array}{lcl} 
1 &\text{si}& \s=1\\
0 &\text{si}& \s\neq1
\end{array}\right.
\end{equation}
L'algèbre $\assl{\QM}{\YB}$ admet une $\NM\times\G$-graduation, 
le {\em poids coloré}, obtenue en convenant que $y_{n,\nu}$ est de degré 
$(n,\nu)$. La première composante, le {\em poids} est héritée de  
$\assl{\QM}{\XB}$. 
La seconde est une $\G$-graduation, qu'on appellera la {\em couleur totale}. 
En convenant que chaque élément de $\YB$ est de degré 1, on définit encore
une autre $\NM$-graduation, la {\em longueur}\footnote{Goncharov utilise le
terme de \og{}profondeur\fg{} (\em{depth}).} 
qui jouera un rôle moins 
important dans la suite. La différence du poids et de la longueur est
le degré partiel en $x_0$ de $\assl{\QM}{\XB}$. 

Le coproduit $\Det$, défini comme morphisme d'algèbres par la condition
\begin{equation}\label{eq_defDet}
\Det(y_{n,\nu})\ \ass\ \underset{\kappa\lambda=\nu}{\sum_{k+l=n}}
y_{k,\kappa}\ot y_{l,\l},
\end{equation}
est clairement cocommutatif, coassociatif et coünifère et fait 
de $(\assl{\QM}{\YB},\pcdot,\Det)$ une bigèbre de Hopf --- l'antipode est
fourni par \ref{sss_U} ci-dessous). Il est 
de plus homogène pour le poids (et même pour le
poids coloré), ce qui permet de le prolonger à $\sernc{\k}{\YB}$. Par contre,
il ne respecte pas la graduation de longueur, mais seulement 
la filtration décroissante associée. 

Il est caractérisé, en tant que morphisme 
d'algèbres homogène pour le poids coloré, par le fait que l'élément 
$$\YC\ \ass\ \sum_{(n,\nu)\in \NM\times\G} y_{n,\nu}$$
de $\sernc{\QM}{\YB}$ est \glk{}, comme le montre un calcul immédiat. 
\sss{Sommes partielles} \label{sss_sommepart}
Les relations du type (\ref{exemple_rel2})
sont déjà vraies pour
les sommes partielles $L^N_{s_1,\ldots,s_r}(z_1,\ldots,z_r)$ définies en 
\ref{sss_divlogLN}.

Fixons $N>0$ et considérons dans $\sernc{\CM}{\YB}$ la série génératrice
\begin{equation}\label{def_LCCN}
\LCC_N\ \ass\ \sum_{r, s_1,\ldots,s_r,\s_1,\ldots,\s_r}
L^N_{s_1,\ldots,s_r}(\s_1,\ldots,\s_r) y_{s_1,\s_1}\ldots y_{s_r,\s_r}
\end{equation}
Ne pas fixer $N$ revient à la voir comme un élément de 
$\sernc{\CM^\NM}{\YB}$, plus précisément de $\sernc{\divlogs}{\YB}$
(cor. \ref{prop_LNdivlog}).

Par test sur les générateurs, on voit que 
pour $\l\in\CM$ et $\nu\in\G$, 
les substitutions $y_{n,\nu}\mapsto \l^ny_{n,\nu}$ et $y_{n,\nu}\mapsto
\nu y_{n,\nu}$ sont des morphismes de cogèbres.\footnote{Pour la première,
il s'agit simplement de l'homogénéité de $\Det$.} 
En appliquant une fois la 
première et $m$ fois la seconde à $\YC$, on obtient donc que l'élément 
$$ \YC_m(\l)\ \ass\ \sum_{n\geq 0,\nu\in\G} \nu^m \l^ny_{n,\nu}$$ de
$\sernc{\CM}{\YB}$ est encore \glk{} pour $\Det$. 

Dans le développement du produit 
$$ P\ \ass\ \YC_{N-1}\left(\frac{1}{N-1}\right)\YC_{N-2}\left(\frac{1}{N-2}\right)
\cdots \YC_{1}\left(\frac{1}{1}\right),$$
les occurences du mot $w\ \ass\ y_{s_1,\s_1}\cdots y_{s_r,\s_r}$ correspondent
aux suites $N>n_1>\cdots>n_r$ d'entiers, le facteur $y_{s_i,\s_i}$ de $w$ 
provenant du facteur $\YC_{n_i}(1/n_i)$ de $P$.
 Le coefficient de $y_{s_i,\s_i}$ dans 
$\YC_{n_i}(1/n_i)$ étant $\s_i^{n_i}n_i^{-s_i}$, l'occurence 
$N>n_1>\cdots>n_r$ de $w$ dans $P$ porte donc le coefficient 
$\s_1^{n_1}n_1^{-s_1}\cdots \s_r^{n_r}n_r^{-s_r}$. Le coefficient total
de $w$ dans $P$ s'obtenant en sommant les coefficients portés par
ces occurences, on voit donc que $P$ est égal à $\LCC_N$. 
Comme chaque $\YC_i(1/i)$ est \glk, on a donc démontré:
\begin{prop}\label{prop_LNdiag}
La série génératrice des sommes partielles $\LCC_N$ est \glke{} 
dans $(\sernc{\CM}{\YB},\Det)$.
\end{prop} 
\sss{Quotient de convergence} 
La formule (\ref{eq_defDet}) exprime en particulier la primitivité de
$y_{1,1}$ dans $(\assl{\QM}{\YB},\Det)$.
L'idéal à droite $I=y_{1,1}\assl{\QM}{\YB}$ est donc
un coïdéal pour $\Det$, le quotient $\assl{\QM}{\YB}_\cv$ hérite du 
coproduit $\Det$ et la projection 
$\pi_\cv:\assl{\QM}{\YB}\to\assl{\QM}{\YB}_\cv$ est un morphisme de 
cogèbres. 
Comme en \ref{sss_convpart}, le projeté $\pi_\cv(\LCC_N)$ est 
\glk{} dans $(\sernc{\CM^\NM}{\YB}_\cv, \Det)$ et on peut passer à la limite:
\begin{prop}\label{prop_Ldiag}
L'élément $\LCC_\cv\ \ass\ \lim\limits_{N\to\infty}\pi_\cv(\LCC_N)$ de 
$\sernc{\CM}{\YB}_\cv$ est \glk{} pour le coproduit $\Det$. 
\end{prop}
\sss{Redressement du coproduit $\Det$} \label{sss_U} 
Malgré la parenté apparente,
au moins au niveau informel, des deux types de relation de mélange,
$(\assl{\QM}{\YB},\pcdot,\Det)$ n'est pas la bigèbre enveloppante de 
$\liel{\QM}{\YB}$, mais elle peut s'y ramener, par un automorphisme 
d'algèbres.

Soit $u_{n,\nu}$ la partie homogène de poids coloré $(n,\nu)$ de 
$\log(\YC)$ (\cf{} \ref{sss_defYDet}) et $\UB$ 
l'ensemble des $u_{n,\nu}, (n,\nu)\in\NM^*\times\G$ (en particulier, ceci 
impose $u_{0,\nu}=0$, pour toute valeur de $\nu$). Comme $\YC$ est \glk{},
son logarithme est primitif; $\Det$ étant homogène pour le poids coloré,
les $u_{n,\nu}$ sont tous primitifs. 

   Par définition, on a 
\begin{equation}
\exp\left(\sum_{n\in\NM,\nu\in\G} u_{n,\nu}\right) = 
\sum_{n\in\NM, \nu\in\G} y_{n,\nu}
\end{equation}
En observant la partie homogène de poids coloré $(n,\nu)$ et de longueur 1
de cette expression, on voit que la partie de longueur 1 de $u_{n,\nu}$ est
$y_{n,\nu}$. On a donc 
\begin{equation}\label{eq_UYlng1}
 u_{s_1,\s_1}u_{s_2,\s_2}\cdots u_{s_r,\s_r} =
y_{s_1,\s_1}y_{s_2,\s_2}\cdots y_{s_r,\s_r} + 
(\text{termes de longueur $>r$})
\end{equation}
En d'autres termes, la restriction à la 
composante homogène de poids $n$ de la 
substitution $y_{n,\nu}\mapsto u_{n,\nu}$ est unipotente. Cela assure 
que cette substitution est un automorphisme d'algèbres de $\assl{\QM}{\YB}$. 
Compte-tenu de la primitivité des $u_{n,\nu}$, on a donc prouvé:
\begin{prop}\label{prop_U}
   La $\QM$-algèbre associative $\assl{\QM}{\YB}$ est librement engendrée par 
$\UB$. La bigèbre $(\assl{\QM}{\YB},\pcdot,\Det)$ est isomorphe à 
$(\assl{\QM}{\UB},\pcdot,\D)$, bigèbre enveloppante universelle de 
$\liel{\QM}{\UB}$. 
\end{prop}

\sss{}On utilisera plus loin la 
dérivation $\d_{u_{s,\s}}$ de $\assl{\QM}{\YB}$
qui envoie $u_{s,\s}$ sur 1 et annule tous les autres $u_{n,\nu}$. Soit $\psi$ 
un 
élément primitif de $(\sernc{\k}{\YB}, \Det)$. D'après la 
proposition ci-dessus, $\psi$ est une série de Lie en les $u_{n,\nu}$. 
On voit facilement que $\d_{u_{s,\s}}(\psi)$ est le coefficient de $u_{s,\s}$
dans $\psi$, car $\d_{u_{s,\s}}$ annule tous les crochets. Compte-tenu de
(\ref{eq_UYlng1}), cela devient:
\begin{prop}\label{prop_deru}
Pour $\psi$ primitif dans $(\sernc{\k}{\YB},\Det)$ et $(s,\s)\in\NM^*\times\G$, on a 
$$\d_{u_{s,\s}}(\psi) = \scal{\psi}{y_{s,\s}}$$
\end{prop}
\sss{Remarques}\label{sss_remqsym}
Dans le cas $\G=\triv$, la bigèbre $(\assl{\QM}{\YB},\cdot,\Det)$ est 
celle des fonctions symétriques non-commutatives \cite{Noncommsym}, 
vue dans la base des fonctions complètes. Sa 
duale graduée est celle des fonctions 
quasi-symétriques $(\qsym{}, \et, \delta)$ 
\cite{MalReut}. Les propositions \ref{prop_LNdiag} et \ref{prop_Ldiag} 
se dualisent en interprétant $L_N$ (resp. $L$), 
après extension par linéarité, comme un morphisme d'algèbres de $\qsym{}$ dans
$\CM$ (resp. $\qsymcv{}\to\CM$, où $\qsymcv{}$ est une sous-algèbre de $\qsym{}$, engendrée
comme $\QM$-espace vectoriel par les éléments de la base duale de 
l'ensemble des mots en $\YB$ 
(les monomiales) 
indicés par $(s_1,\ldots,s_r)$ avec $s_1\neq 1$; elle est duale du quotient 
$\assl{\QM}{\YB}_\cv$. Pour effectuer des calculs
explicites, il est souvent préférable de travailler avec $\qsym{}$.  

Pour généraliser  à $\G$ 
quelconque, M. Bigotte \cite{Big} a introduit l'algèbre des 
\og{}fonctions quasi-symétriques colorées\fg{}. On laisse au lecteur 
le soin de vérifier qu'elle est duale de la cogèbre 
$(\assl{\QM}{\YB_\G}, \Det)$. 

La proposition \ref{prop_U} est une généralisation directe du
cas particulier $\G=\triv$, classique pour les fonctions symétriques 
non-commutatives. On a des descriptions plus fines de ce changement de base,
ainsi que de nombreuses variantes \cite{Noncommsym,MalReut}.  
On pourrait encore généraliser en rempla\c cant 
$\NM^*\times\G$ par un magma associatif et commutatif $(M,+)$ quelconque, 
dans lequel tout élément aurait un nombre fini de décompositions en sommes
de deux éléments. La convention (\ref{conv_y0}) se formulerait en 
adjoignant un élément neutre à $M$. 

\sss{}Nous aurons dans la suite besoin d'un calcul qui se fait 
habituellement dans $\qsym{}$. Appliqué à l'élément \glk{} $\LCC_\cv$, 
il redonne l'exemple \ref{exemple_rel2} de l'introduction (voir 
\ref{subsec_conv} pour la notation). 

\begin{prop}\label{prop_etpart}
   Pour tous $v\in\sernc{\k}{\YB}, y_{s,\s}\in \YB$ et $y_{t,\tau}\in \YB$,
on a:
$$ \scal{\Det v}{y_{s,\s}\ot y_{t,\tau}} = \scal{v}{y_{s+t,\s\tau} + 
y_{s,\s}y_{t,\tau} + y_{t,\tau}y_{s,\s}} $$
Pour $v$ primitif pour $\Det$,  $s\neq0$ et $t\neq0$, on a donc 
\begin{equation}\label{eq_etpart}
 \scal{v}{y_{s+t,\s\tau} + y_{s,\s}y_{t,\tau} + y_{t,\tau}y_{s,\s}} = 0
\end{equation}
\end{prop}
\begin{preuve} Comme $\Det$ respecte la filtration de longueur, 
seuls les termes de longueur 1 et 2 de $v$ interviennent. 
L'identification est alors immédiate. 
La conséquence (\ref{eq_etpart}) est évidente.
\end{preuve}
\subsection{Régularisations} Les séries $\ItC_\cv$ et $\LCC_\cv$ sont
peu pratiques à utiliser, car les structures de bigèbre ne persistent pas
dans les quotients de convergence. Il est donc bénéfique
de les relever à des éléments diagonaux des deux bigèbres concernées.  
D'un point de vue combinatoire, ce processus n'est pas canonique,
mais on a un bon contrôle de l'ambiguïté,
permettant de fixer arbitrairement des relevés $\ItC$ et $\LCC$.  

On peut considérer les coefficients apparaissant dans $\ItC$ et $\LCC$ comme
des valeurs régularisées des intégrales itérées et des sommes itérées
préservant chacune la relation de mélange adéquate.
Il est bien connu que la formule de Kontsevitch ne survit pas au processus. 
Dans notre formalisme, cela s'écrit $\LCC\neq\qs\pi_Y(\ItC)$. 
 
Les deux types de 
développements asymptotiques de \ref{subsec_divlog}
fournissent d'autres relevés, à coefficients polynômiaux. 
Les propriétés de comparaison de \ref{sss_comp} se traduisent,
grâce au contrôle de l'ambiguïté des relèvements, par la relation de 
régularisation, qui exprime $\LCC$ en fonction de $\ItC$. 
On évoque ensuite ses conséquences et son historique.      

\sss{Cadre général}
Soient $\gG$ une $\QM$-algèbre de Lie graduée, $\gG_1$ et $\gG_2$ deux 
sous-algèbres de Lie de $\gG$, telles que $\gG=\gG_1\oplus\gG_2$, en tant 
qu'espaces vectoriels.  

Le théorème de Poincaré-Birkhoff-Witt montre que l'application 
\begin{equation}
\text{PBW} : \begin{array}{ccc}
\Uec{}{\gG_1}\ot\Uec{}{\gG_2}&\longrightarrow&\Uec{}{\gG}\\
a\ot b&\longmapsto&ab
\end{array}
\end{equation}
est un isomorphisme de cogèbres. Notons $\eps$ la coünité de $\Uec{}{\gG_1}$.
L'application  $\eps\ot\id : \Uec{}{\gG_1}\ot\Uec{}{\gG_2}\to\Uec{}{\gG_2}$ 
est un morphisme de cogèbres, dont le noyau 
$\gG_1\Uec{}{\gG_1}\ot\Uec{}{\gG_2}$ a pour image 
$\gG_1\Uec{}{\gG}$ par PBW.
La cogèbre-quotient $\Uec{}{\gG}/\gG_1 \Uec{}{\gG}$ est
donc isomorphe à $\Uec{}{\gG_2}$ et on a un diagramme commutatif de 
cogèbres:
\begin{equation}
\xymatrix{
\Uec{}{\gG_1}\ot\Uec{}{\gG_2}\ar@<-1ex>[d]_-{\eps\ot\id}\ar[r]^-{\text{PBW}}&\Uec{}{\gG}\ar[d]^{\pi} \\
\Uec{}{\gG_2}\ar@<-1ex>[u]_-{1\ot\id}\ar[r]^{~}&\Uec{}{\gG}/\gG_1 \Uec{}{\gG}
}
\end{equation}
Les deux flèches horizontales de ce diagramme étant des isomorphismes de 
cogèbres, on déduit de $1\ot\id$ une section de la projection canonique $\pi$,
avec comme conséquence:
\begin{prop}\label{prop_regulgen}
Tout élément \glk{} de $\mglk{}(\Uec{\k}{\gG}/\gG_1\Uec{\k}{\gG})$ 
est l'image par $\pi$ d'un élément de $\exp(\k\gG)$. 
Deux tels éléments se déduisent l'un de l'autre par multiplication à gauche
par un élément de $\exp(\k\gG_1)$. 
\end{prop}
On utilisera sans autre commentaire la variante pour les multiplications
à droite. 
\sss{Application à $\ItC_\cv$ et $\LCC_\cv$}
\vdeux{%
D'après le théorème d'élimination de Lazard \cite{Reut}\fnote{Préciser}, 
$\liel{}{\UB}$ est somme directe de $\QM u_{1,1}$ et de l'idéal engendré par 
les $\ad u_{1,1}^k(u_{n,\nu})$ pour $k>0$ et $(n,\nu)$ quelconque.}
{La proposition \ref{prop_regulgen} s'applique notamment lorsque $\gG$ est
une algèbre de Lie libre sur un alphabet $\AB$ et 
$\gG_1=\k a,$ pour $a\in\AB$: on prend pour $\gG_2$ le premier terme non 
trivial de la filtration décroissante associée au degré partiel en 
$\AB\moins{a}$.}    
 
La bigèbre enveloppante de $\liel{}{\UB}$ est $(\assl{\QM}{\YB},\pcdot,\Det)$
(prop. \ref{prop_U}) et on a $y_{1,1}=u_{1,1}\in\UB$. 
Comme les coefficients de $y_{1,1}$ 
s'additionnent lorsqu'on multiplie deux séries \glkes{}, on en déduit:
\begin{cor}\label{cor_regulDet}
   Toute série $\Phi_\cv$ \glke{} de $(\assl{\QM}{\YB}_\cv,\Det)$ 
   est l'image par $\pi_\cv$ 
d'un élément \glk{} de $(\assl{\QM}{\YB}, \Det)$. 
   Deux tels éléments $\Phi_1$ et $\Phi_2$ sont liés par 
   $$ \Phi_2 = \exp((\l_2-\l_1)y_{1,1})\Phi_1, $$
   où $\l_i$ est le coefficient de $y_{1,1}$ dans $\Phi_i$, 
   pour $i\in\{1,2\}$. 
\end{cor}
\vdeux{Pour $\gG=\liel{}{\XB}$, on a 
$\gG=\QM x_1\oplus [\gG,\gG]\oplus \QM x_0$. En appliquant deux fois la 
proposition \ref{prop_regulgen}, }
{De même, on peut décomposer $\liel{}{\XB}$ en 
$\QM x_1\oplus\hG\oplus\QM x_0$ et appliquer deux fois la proposition 
\ref{prop_regulgen}; }
il vient d'abord une variante pour 
$(\assl{\QM}{\YB},\D)$ du corollaire \ref{cor_regulDet}, puis: 
\begin{cor}\label{cor_regulD}
   Toute série $\Phi_\cv$ \glke{} de $(\assl{\QM}{\XB}_\cv,\D)$ 
   est l'image par $\pi_\cv$ d'un élément \glk{} de $(\assl{\QM}{\XB}, \D)$. 
   Deux tels éléments $\Phi_1$ et $\Phi_2$ sont liés par 
   $$ \Phi_2 = \exp((\l_2-\l_1)x_1)\Phi_1\exp((\mu_2-\mu_1)x_0), $$
   où $\l_i$ et $\mu_i$ sont respectivement les coefficients de $x_1$ et $x_0$ 
   dans $\Phi_i$, pour $i\in\{1,2\}$. 
\end{cor}    
Ces corollaires peuvent également se déduire \cite{PetMinh,thesejoe}  
d'un processus plus fin de décomposition des exponentielles de Lie, 
bien plus adapté pour les calculs informatiques que 
les développements en combinaisons linéaires de mots. 

\sss{Notations}
On désignera par $\ItC$ l'unique série \glke{} de $(\sernc{\CM}{\XB},\D)$ 
telle que $\pi_\cv(\ItC)=\ItC_\cv$ 
et dans laquelle les coefficients de $x_0$ et $x_1$ sont nuls. De même, 
l'unique élément diagonal de $(\sernc{\CM}{\YB},\Det)$ se projetant par 
$\pi_\cv$ sur $\LCC_\cv$ et dans laquelle le coefficient de $y_{1,1}$ est nul
sera noté $\LCC$. L'existence et l'unicité de ces deux séries 
sont garanties par les corollaires \ref{cor_regulD} et \ref{cor_regulDet}. 

La série $\pi_Y(\ItC)$ est caractérisée par le fait qu'elle est 
\glke{} dans $(\sernc{\CM}{\YB},\D)$, antécédent de $\ItC_\cv$ par $\pi_\cv$
et que le coefficient de $y_{1,1}$ y est nul.

La suite de cette partie est consacrée à l'obtention d'une formule explicite de
comparaison entre $\ItC$ et $\LCC$.

On adoptera la convention suivante: soient $\k_1$ et $\k_2$ deux 
$\CM$-anneaux,\fnote{Dans le par. de conventions ?}
$f : \k_1\to\k_2$ une application $\CM$-linéaire et 
$\Phi\in\sernc{\k_1}{\XB}$. On note $\ov{f}(\Phi)$ l'élément de  
$\sernc{\k_2}{\XB}$ obtenu en appliquant $f$ aux coefficients de $\Phi$. 
L'application $\ov f$ ainsi définie est un morphisme de 
$\sernc{\CM}{\XB}$-bimodules. Si $f$ est un morphisme d'algèbres et
$\Phi$ est \glk{} pour $\D$, alors $\ov{f}(\Phi)$ l'est aussi. De même
pour les $\YB$ et $\Det$.

\sss{Interprétation pour les sommes itérées}\label{sss_regulaset}
Les coefficients de $\LCC_N$ sont des suites de $\divlogs$ 
(prop. \ref{prop_LNdivlog}). 
Comme $\As$ est un morphisme d'algèbres, $\LCC_t\ \ass\ \ov{\As}(\LCC_N)$
est \glk{} dans $(\sernc{\CM[t]}{\YB},\Det)$ ; il se projette de plus par
$\pi_\cv$ sur $\LCC_\cv$, si l'on considère $\sernc{\CM}{\YB}$ comme inclus 
dans $\sernc{\CM[t]}{\YB}$. 
Le coefficient de $y_{1,1}$ dans $\LCC_t$ étant $\As((H_N))$, égal à 
$\g+t$, le corollaire \ref{cor_regulDet} 
appliqué avec $\k=\CM[t]$ détermine $\LCC_t$:
\begin{equation}\label{eq_LCCt}
   \LCC_t = \exp((\g+t)y_{1,1})\LCC
\end{equation}

\sss{Justifications de \ref{sss_comp}}\label{sss_aposteriori}
Par $1^k$, on entend la séquence formée de $k$ fois le nombre $1$. 
Par définition, $H_{k,N}\ \ass\ L^N_{1^k}(1^k)$ est le
coefficient de $y_1^k$ dans $\LCC_N$. Le polynôme $\As((H_{k,N}))$ est
donc le coefficient de $y_1^k$ dans $\LCC_t$. La formule (\ref{eq_LCCt}) 
montre qu'il est de degré $k$ en $t$.  

D'après la formule de Kontsevitch, $L_{1^k}(1^{k-1},z)$, valeur en $z$  
de la fonction associée à $H_{k,N}$, n'est autre que le coefficient de $x_1^k$
dans l'élement \glk{} $\ItC_{[0,z]}$. 
D'après le lemme ci-dessous, c'est donc $(L_1(z))^k/k!$, car le 
coefficient 
de $x_1$ dans $\ItC_{[0,z]}$ est $L_1(z)$.     
\begin{lemme}\label{lemme_coeffx1}
Soient $\k$ un $\QM$-anneau, $\Phi$ un élément \glk{} de 
$(\sernc{\k}{\XB},\D)$ et 
$\l$ le coefficient de $x_1$ dans $\Phi$. Le coefficient de 
$x_1^k$ dans $\Phi$ est $\l^k/k!$.
\end{lemme}     
\begin{preuve}
Munie du coproduit $\D$ pour lequel $t$ est primitif, $\serc{\QM}{t}$ 
est la bigèbre enveloppante universelle complétée 
de l'algèbre de Lie commutative de base $\{t\}$. 
La substitution $\pi_t: x_1\mapsto t, x_\alpha\mapsto 0$ pour 
$\alpha\neq 1$ est clairement un morphisme de cogèbres de 
$(\sernc{\k}{\XB},\D)$ dans $(\serc{\k}{t},\D)$. L'image de $\Phi$ par 
$\pi_t$ est donc \glke{}, \ie{} de la forme $\exp(\mu t)$. 
La comparaison des termes de degré 1 donne $\l=\mu$. 
\end{preuve}

\sss{Interprétation asymptotique pour les intégrales itérées}
\label{sss_asympit}
Les coefficients de $\LCC_z$ et de $\ItC_{[0,z]}$ sont les fonctions 
$z\mapsto L_{s_1,\ldots,s_r}(\s_1,\ldots,z\s_r)$, 
associées comme en \ref{sss_serfonc} aux suites 
$(L^N_{s_1,\ldots,s_r}(\s_1,\ldots,\s_r))_{N>0}$, qui sont dans $\divlogs$.
D'après le  
corollaire \ref{cor_comp}, dont on vient de compléter la démonstration,
ce sont donc des éléments de $\divlogf$.

On peut donc poser $\ItC_t\ \ass\ \ov\As(\ItC_z)$ ; c'est un élément 
diagonal de
$(\sernc{\CM[t]}{\YB}, \D)$ qui se
projette par $\pi_\cv$ sur $\ItC_\cv$. Le coefficient de $x_1$ dans 
$\ItC_t$ vaut cette fois $-t$, car c'est $\As(L_1)$.  
Le corollaire \ref{cor_regulD} nous donne donc:
\begin{equation}\label{eq_ItCt}
\ItC_t=\exp(-x_1t)\pi_Y(\ItC)
\end{equation}
On voit au passage que $\pi_Y(\ItC)$ s'obtient ---  ainsi donc que 
les valeurs régularisées des intégrales itérées de $0$ à $z$  --- 
en prenant le terme 
constant des développements asympotiques polynomiaux en $\log(1-z)$ 
des intégrales itérées de $0$ à $z$ au voisinage de $z=1$. On retrouve 
bien la \og{}régularisation canonique\fg{} de Goncharov, et 
la description de l'image d'un chemin aboutissant à un point-base 
tangentiel dans la comparaison Betti-de Rham \cite[15.52]{DelPi1}.
\sss{Comparaison des deux régularisations}\label{sss_compregul}
On a $\LCC_z=\qs\ItC_{[0,z]}$ et $\LCC_z$ se déduit de 
$\LCC_N$ par l'application $\CM$-linéaire 
$\divlogs\to\divlog_D$ de \ref{sss_comp}.

D'après les définitions de $\ItC_t$ et $\LCC_t$ et le corollaire \ref{cor_comp},
on a $\ov{\comp}(\LCC_t)=\qs\ItC_t$. Comme $\ov\comp$ est un endomorphisme de
$\sernc{\CM}{\YB}$-modules à droite, 
on déduit de (\ref{eq_LCCt}) et (\ref{eq_ItCt}) 
que l'on a $\LCC=\LCC_\corr \qs\pi_Z(\ItC)$, avec 
$\LCC_\corr\ \ass\ \ov\comp(\exp((\g+2t)y_1))$. La série $\LCC_\corr$ est
 un produit d'éléments de $\sernc{\CM}{\YB}$ ; elle est donc à coefficients 
dans $\CM$. On a ainsi prouvé: 
\begin{prop}\label{prop_compregul}
 Il existe $S\in\serc{\CM}{y_{1,1}}$ telle que
que $$\LCC=S\cdot\qs\pi_Y(\ItC)$$ \end{prop}
\sss{Relation de régularisation} \label{sss_relregul}
La proposition \ref{prop_compregul} suffit à déterminer de manière purement
algébrique les coefficients de $S$. 

\begin{prop}\label{prop_fixecorr}
   Soit $\k$ un $\QM$-anneau, $\Phi$ un élément \glk{} de 
$(\sernc{\k}{\XB},\D)$, et une série $S$ de $\serc{\k}{t}$ telle
que l'élément $\Phi_\et\ \ass\ S(y_{1,1})\cdot\qs\pi_Y(\ItC) $
de $\sernc{\k}{\YB}$ soit \glk{} pour $\Det$. On a 
$$  S = \exp\left(\sum_{n\geq 1} \frac{(-1)^{n-1}}{n}\scal{\Phi}{x_0^{n-1}x_1}
t^n\right) $$
\end{prop}
\begin{preuve}
L'application $\pi_t$ du lemme \ref{lemme_coeffx1} 
se factorise par $\pi_Y$. Le quotient, noté encore 
$\pi_t$, est la substitution $y_{1,1}\mapsto t, y_{n,\nu}\mapsto 0$ pour 
$(n,\nu)\neq (1,1)$. 
Par le lemme 
\ref{lemme_coeffx1}, on a 
$\pi_t\qs\pi_Y(\Phi)=1$, d'où $\pi_t(\Phi_\et)=S$.  

D'après la proposition \ref{prop_U}, on peut écrire $\Phi$ comme 
exponentielle d'une série de Lie en les $u_{n,\nu}$. Pour $n\neq0$, notons 
$\alpha_{n,\nu}$ le coefficient de $u_{n,\nu}$
dans celle-ci. 

Comme $\serc{\CM}{t}$ est commutative, tout crochet
des $u_{n,\nu}$ est dans le noyau de $\pi_t$, d'où 
\begin{equation}\label{eq_LCCun}
\log\pi_t(\Phi_\et)= \sum_{n\in\NM^*, \nu\in\G}\alpha_{n,\nu}\pi_t(u_{n,\nu})
\end{equation}
On a $\pi_t(u_{n,\nu})=0$ si $\nu\neq 1$, car 
$u_{n,\nu}$ est homogène de couleur totale $\nu$. On en déduit
$$ \sum_{n>0} \pi_t(u_{n,1}) = 
\pi_t\left(\sum_{n\in\NM, \nu\in\Gamma} u_{n,\nu}\right)\ \egdef\ 
\pi_t\log\left(\sum_{n\in\NM, \nu\in\Gamma}y_{n,\nu}\right) = \log(1+t),$$
d'où l'on tire $\pi_t(u_{n,1})=(-1)^{n-1}t^n/n$ par homogénéité de $\pi_t$.
Reportons dans (\ref{eq_LCCun}):
\begin{equation}\label{eq_LCCan}
  \log\pi_t(\Phi_\et)=\sum_{n>0}\frac{(-1)^{n-1}}{n}\alpha_{n,1} t^n
\end{equation}
Il reste à déterminer les $\alpha_{n,1}$. À des termes de longueur 
au moins 2 près, on a $u_{n,\nu}=y_{n,\nu}$ (\cf{} \ref{sss_U}) et donc 
$\Phi_\et=\sum_{n,\nu} \alpha_{n,\nu}y_{n,\nu}$, car une série de Lie n'a pas de 
terme constant. On voit donc que 
$\alpha_{n,1}$ est le coefficient de $y_{n,1}$ dans $\Phi_\et$, \ie{}
$\scal{\Phi}{x_0^{n-1}x_1}$. 
\end{preuve}
Le coefficient de $x_0^{n-1}x_1$ dans $\ItC$ vaut $\z(n)$ si $n\geq 2$ et 
0 si $n=1$. 
\begin{cor}[Relation de régularisation] \label{prop_relregul} 
Les séries $\LCC$ et $\ItC$ vérifient
$$\LCC= \exp\left(\sum_{n\geq2}\frac{(-1)^{n-1}\z(n)}{n}y_{1,1}^n\right)\qs\pi_Y(\ItC)$$
\end{cor}
\sss{Conséquences} 
Les relations DMR (doubles mélanges et régularisation) permettent de démontrer
l'egalité $\z(2,1)=\z(3)$, qui était connue d'Euler, et ne peut 
se déduire des seuls doubles mélanges, relations de 
poids minimum 4.  

Les relations DMR impliquent\footnote{Voir \cite{thesejoe} cor. III.4.21 p. 97
pour le détail, dans 
le cas des polyzêtas.} plus généralement 
la relation d'Hoffman \cite{Hoffalg}, qui se généralise 
aux racines de l'unité \cite{Big} et sur laquelle les calculs de 
\cite{PetMinh, Big} se fondent, en plus des doubles mélanges: elle se ramène
à l'absence de terme en $y_{1,1}$ dans le facteur correctif. 
On sait peu sur la réciproque, si ce n'est qu'elle est vraie pour les 
polyzêtas aussi loin que 
l'on puisse calculer par ordinateur avec les deux systèmes, 
\ie{} jusqu'en poids 16.
\sss{Remarques} Le calcul de la proposition \ref{prop_fixecorr}
se ramène à la variante de Waring des formules de Newton, 
modulo l'inclusion que nous n'avons pas détaillée  
de l'algèbre des fonctions symétriques dans celle des quasi-symétriques, 
elle-même incluse dans celle des fonctions quasi-symétriques colorées (\cf{}
\ref{sss_remqsym}). 

	L'attribution de la relation de régularisation est difficile,
même en faisant abstraction des différences de langage. Elle apparaît chez 
Goncharov\footnote{Ceci nous avait échappé à l'époque de la rédaction de  
\cite{thesejoe}.} \cite{Gonch98, Gonch2001}, sous une forme duale de la proposition 
\ref{prop_compregul} ; chez Écalle \cite{Ecalle}, par des arguments très 
différents, modulo la correspondance 
entre ses moules entiers et nos séries génératrices, sous la forme 
\ref{prop_compregul} ; chez Boutet de Monvel \cite{Boutet}, avec la définition
et le calcul complet de $\comp$. Plus  
récemment, on trouve dans \cite{Kaneko} des arguments proches,  
cette fois attribués à Zagier. 
   
\Subsection{Relations spécifiques aux racines de l'unité} 
\sss{} 
Les relations de distribution du polylogarithme classique 
se généralisent directement aux polylogarithmes multiples \cite{Gonch2001}:
\begin{prop}\label{prop_distclass}
Soit $d$ un diviseur de l'entier $n$. Si $(s_1,\s_1)\neq(1,1)$, on a 
\begin{equation}\label{eq_distrib}
 \sum_{\t_1^d=\s_1,\t_2^d=\s_2,\ldots,\t_r^d=\s_r}
L_{s_1,\ldots,s_r}(\t_1,\ldots,\t_r)
= d^{r-(s_1+\cdots+s_r)}L_{s_1,\ldots,s_r}(\s_1,\ldots,\s_r)
\end{equation}
\end{prop}
Il s'agit de traduire cela en termes de séries génératrices non-commutatives.
\sss{Fonctorialités en $\G$}\label{sss_proj}
Tout morphisme de groupes $\phi:\G\to\G'$ donne lieu à deux substitutions 
$\phi_*:\assl{\QM}{\XB_\G}\to\assl{\QM}{\XB_{\G'}}$ et 
$\phi^*:\assl{\QM}{\XB_{\G'}}\to\assl{\QM}{\XB_\G}$, définies par:
\begin{equation}
\phi^*(x_\s) = \left\{
\begin{array}{lcl}
x_0 &\text{si}& \s=0 \\
\sum\limits_{\tau\in\phi^{-1}(\s)}x_\tau &\text{si}& \s\in\G \\
\end{array}\right. \qtext{et}
\phi_*(x_\s) = \left\{
\begin{array}{lcl}
d x_0 &\text{si}& \s=0\\
x_{\phi(s)} &\text{si}& \s\neq 0, 
\end{array}\right. 
\end{equation}
où $d$ désigne l'ordre de $\ker\phi$. 

Ces applications préservent les structures qu'on a précédemment définies. 
Ce sont évidemment des morphismes de cogèbres pour $\D$. Elles stabilisent  
les $\assl{\QM}{\YB}$, commutent à $\pi_Y$. Par test sur les 
$y_{n,\nu}$, on voit que ce sont des morphismes de cogèbres pour $\Det$. 
L'image directe $\phi_*$ commute à 
$\ps$ et $\qs$; si $\phi$ est injective, l'image réciproque $\phi^*$
stabilise $\assl{\QM}{\YB}_\cv$ et commute à $\ps,\qs$ et $\pi_\cv$. 

\sss{}\label{sss_distrib}
Soit $\G$ un groupe commutatif fini. 
Pour tout diviseur $d$ de l'ordre de $\G$, notons $p^d$ l'application 
$\s\mapsto\s^d$, $\G^d$ son image et $i_d$ l'inclusion de $\G^d$ dans $\G$.

Pour $\G\subset\CM^*$, on laisse au lecteur le soin de
vérifier que la proposition \ref{prop_distclass} 
se traduit par l'égalité, dans $\sernc{\CM}{\YB_{\G^d}}$, de 
$i_d^*(\LCC_\cv)$ et $\pi_\cv p^d_*(\LCC)$. La même égalité 
vaut pour $\ItC$, car $i_d^*$ et $p^d_*$ commutent à
$\ps$ et $\qs$. Les séries 
$i_d^*(\ItC)$ et $p^d_*(\ItC)$, toutes 
deux \glkes{}, ont donc même image, $i_d^*(\LCC_\cv)$,
par $\pi_\cv$. D'après
le corollaire \ref{cor_regulD},  
il suffit de connaître leurs termes en $x_0$ et $x_1$ 
pour les comparer ; on en déduit facilement 
la relation de distribution régularisée: 
\begin{prop}
Pour tout sous-groupe multiplicatif fini $\G$ de $\CM^*$ et tout
diviseur $d$ de l'ordre de $\Gamma$, on a:
$$  p^d_*(\ItC) =
\exp\left(\sum_{\s^{n/d}=1, \s\neq1}\Li_1(\s)x_1\right)i_d^*(\ItC) $$
\end{prop}

\sss{Relations de poids un} \label{sss_poids1}
Un calcul direct permet d'obtenir, pour $k$ strictement compris entre
$0$ et $n/2$:
\begin{equation}
   L_1(\exp(2ki\pi/n))-L_1(\exp(-2ki\pi/n))) = (n-2k)i\pi
\end{equation}
Lorsque $\nu$ varie dans $\mub_n(\CM)$, 
les $L_1(\nu)-L_1(\nu^{-1})$ sont donc tous colinéaires sur $\QM$, 
avec des coefficients explicites, mais qui  
ne sont pas invariants par automorphismes de $\mub_n(\CM)$, à
l'exception de la conjugaison complexe.%
\footnote{Si $n\not\in\{1,2,3,4,6\}$, auquel cas $\mub_n(\CM)$ a plus de deux 
automorphismes, on voit donc que 
les relations de poids un ne découlent pas des autres, qui sont 
invariantes grâce aux propriétés énumérées en \ref{sss_proj}. 
} 
  
D'après \cite{Del01}, ceci traduit la dépendance non fonctorielle 
de la correspondance Betti-de Rham du choix d'un plongement du corps
cyclotomique dans $\CM$. Par contre, la dépendance par rapport au choix 
d'une clôture algébrique de $\RM$ est fonctorielle, 
donnant lieu à l'invariance par conjugaison, qui a 
d'intéressantes conséquences (voir \ref{subsec_conseq}).

\section{Étude formelle des relations DMRD}\label{sec_formel}
 \begin{defi}\label{def_DM}
   Pour tout $\QM$-anneau $\k$ et tout groupe commutatif fini
   $\Gamma$, on note $\DMR(\Gamma)(\k)$ l'ensemble des {\'e}l{\'e}ments $\Phi$
   de $\sernc{\k}{\XB_{\G}}$ tels que:
\begin{eqnarray}
  (\Phi|1) =  1 &\text{\rm et}& (\Phi|x_0) = (\Phi|x_1) = 0
  \label{eq_Phipoids1}\\ 
  \D\Phi\ =\ \Phi\whot_{\k}\Phi &\text{\rm et}&  \Det\Phi_\et\ =\
  \Phi_\et\whot_{\k}\Phi_\et, \\
\qtext{o{\`u} l'on pose} \Phi_\et\ \ass\ \Phi_\corr\cdot\qs\pi_Y(\Phi)&\text{\rm et}&\Phi_\corr\ \ass
 \   \exp\left(\sum_{n\geq2}
{\frac{(-1)^{n-1}}{n}(\pi_Y(\Phi)|y_{n,1})y_{1,1}^n}\right)\nonumber
\end{eqnarray}

Pour tout plongement $\iota$ de $\G$ dans $\CM^*$, on note
$\DMRP(\iota)(\k)$ l'ensemble des éléments $\Phi$ de $\DMR(\G)(\k)$ 
vérifiant les relations de poids un:
\begin{equation}
   \scal{\Phi}{x_{\xi^k}-x_{\xi^{-k}}} = 
\frac{n-2k}{n-2}\scal{\Phi}{x_\xi-x_{\xi^{-1}}},
\end{equation} 
avec $\xi\ \ass\ \iota^{-1}(exp(2i\pi/n))$ et $k$ compris entre 1 et n/2. 

Soit $n$ l'ordre de $\G$. On note respectivement 
 $\PDMRD(\G)(\k)$ et $\DMRD(\iota)(\k)$ 
l'ensemble des éléments $\Phi$ de $\DMR(\G)(\k)$ et de $\DMRP(\iota)(\k)$
vérifiant pour tout diviseur $d$ de $n$ la relation de distribution:
\begin{equation}\label{dist_DMRD}
  p^d_*(\Phi) =
\exp\left(\sum_{\s^{n/d}=1}(\Phi|x_\s)x_1\right)i_d^*(\Phi) 
\end{equation}
Lorsqu'il n'y a pas d'ambiguïté sur $\G$ ou $\iota$, 
on le supprimera de la notation.
\end{defi}

Le contenu de la section \ref{sec_descr} se résume donc ainsi:
pour tout $\iota$, la série $\ItC(\iota(\G))$, 
définie à partir des valeurs prises par les
polylogarithmes multiples sur $\iota(\G)$, est un élément 
de $\DMRD(\iota)(\CM)$. Bien sûr, la donnée de $\iota$ n'est possible que
si $G$ est cyclique, et équivaut à celle de $\xi$. Nombre de propriétés
se formuleront néanmoins pour $\G$ quelconque. 

\begin{defi}\label{def_DMl}
   Soient $\k$, $\iota: \G\to\CM^*$, $\xi$ et $n$ comme précédemment. 
On note resppectivement $\DMRP_\l(\k)$ et $\DMRD_\l(\k)$
   l'ensemble des éléments $\Phi$ de $\DMRP(\k)$ et $\DMRD_\l(\k)$ 
tels que $(\Phi|\alpha_\iota) = \l$, avec 
$$\alpha_\iota\ \ass\ \left\{\begin{array}{lcl}
	x_0x_1 = y_{2,1} &\text{\rm si}& n\in\{1,2\} \\
	\frac{2n}{n-2}(x_\xi-x_{\xi^{-1}})&\text{\rm si}& n\geq3
\end{array}\right.
$$	 
\end{defi}
Plus géométriquement, $\Phi\mapsto\scal{\Phi}{\alpha_\iota}$ 
est un morphisme de schémas $\DMRP\to\AM^1$ 
et  $\DMRP_\l(\k)$ est l'ensemble des $\k$-points de la fibre 
de $\DMRP$ au-dessus de $\l : \spec\k\to\AM^1$. 

Pour $n\leq 2$ (auquel cas $\iota$ est unique), 
le coefficient de $x_0x_1$ dans $\ItC$ est 
par définition $\z(2)$, et donc cette série appartient à 
$\DMRD_{\z(2)}(\CM)$. 

Pour $n\geq3$, le coefficient de $\iota(x_{\xi^k})$ dans $\ItC$ 
vaut $L_1(\exp(2ki\pi/n))$. D'après \ref{sss_poids1}, on a donc dans ce cas 
$\ItC\in\DMRD_{2i\pi}(\CM)$,

\def\mon{\text{\sf M}}
\def\MTG{\text{\sf E}}
\def\MTD{\text{\sf F}}
\def\MTS{\text{\sf S}}
\def\homp#1#2{\hom{\Pi(\k)}{#1}{#2}} 
	
\subsection{Le groupe $\MT$}\label{subsec_MT}
On décrit ici, suivant \vdeux{}{la présentation de }\cite{Del01},
le groupe par lequel se factorise l'action du groupe motivique $\UC_\gr$ 
sur $\pi_1(\Punrac{n})$.
\vdeux{}{Différentes variantes de ces formules apparaissent dans la
littérature, provenant de la composition des  automorphismes 
(extérieurs s'il y a des problèmes d'invariance de points-base) 
\og{}spéciaux\fg{} de $\pi_1(\Punrac{n})$, équivariants pour certaines 
symétries. On ne fait ici que les considérer dans un domaine maximal 
de définition.}

\sss{Premières définitions} 
Soit $\Pi(\k)$ la catégorie dont les objets sont $O_\s$, pour $\s\in\GZ$ et 
avec $\homp{O_\s}{O_\t}=\sernc{\k}{\XB}$, 
la composition des flèches 
étant la multiplication des séries. Un élément $G$ de $\sernc{\k}{\XB}$ sera
noté $G_{\t,\s}$ lorsqu'il est vu dans $\homp{O_\s}{O_\t}$. 

Pour $\s\in\G$, notons $t_\s$ la substitution $x_\nu\mapsto x_{\s\nu}$ de 
$\sernc{\k}{\XB}$. Soit $C_\s$ l'endofoncteur de $\CC(\k)$ 
donné sur les objets par $O_\nu\mapsto O_{\s\nu}$ et 
sur les flèches par $t_\s$.   

\vdeux{}{Dans le cas $\G=\mub_n$, 
la catégorie $\Pi(\k)$ est la $\k$-algèbre du groupoïde fondamental de 
$\Punrac{n}$, en réalisation \og{}graduée\fg{}. Les $\k$-points de ce dernier 
forment la sous-catégorie dont les flèches sont les exponentielles 
de Lie. Les $C_\s$ proviennent de l'action de $\mu_n$ sur $\Punrac{n}$ par
rotations.}

Soit $\MTG(\k)$ l'ensemble des endofoncteurs
de $\Pi(\k)$ agissant trivialement sur les objets, de manière $\k$-linéaire 
et continue sur les flèches, qui fixent les $(x_\s)_{\s,\s}$ et qui commutent
aux $C_\s$. À $\FC\in\MTG(\k)$, faisons correspondre l'élément 
$\Phi(\FC)=\FC(1_{1,0})$. 
On définit ainsi une application de $\MTG(\k)$ dans $\sernc{\k}{\XB}$. 

\sss{}
Fixons $\FC\in\MTG(\k)$. Comme les flèches $1_{0,1}$ et $1_{1,0}$ sont 
inverses l'une de l'autre, $\Phi(\FC)$, on voit que $\FC(1_{0,1})$ est
l'inverse de $\Phi(\FC)$. Pour $\s\in\G$, de $C_\s(1_{1,0})=1_{\s,0}$ découle
$\FC(1_{\s,0})=t_\s(\Phi(\FC))$. 

Soit $\kappa_G$ la substitution  
$x_\s\mapsto t_\s(\Phi(\FC))x_\s [t_\s(\Phi(\FC))]^{-1}$, $x_0\mapsto x_0$. 
L'action de $\FC$ sur $\endo{\Pi(\k)}{x_{0,0}}$ est un morphisme continu de 
$\k$-algèbres. De $(x_\s)_{0,0}=1_{\s,0}(x_\s)_{\s,\s}1_{0,\s}$ 
découle immédiatement que c'est $\kappa_{\Phi(\FC)}$. 

Pour $G,H\in\sernc{\k}{\XB}$, posons, si $G$ est inversible,  
\begin{equation}
\label{eq_defMT}
G\pmt H\ \ass\ G \kappa_G(H) 
\end{equation}
En écrivant $H_{1,0}=1_{1,0}H_{0,0}$, on obtient $\FC(H_{1,0})=\Phi(\FC)\pmt H$. 
Il est clair que l'application $G$ est une bijection de $\MTG(\k)$ 
sur l'ensemble des séries inversibles de $\sernc{\k}{\XB}$. Transportée
par cette bijection, la composition de $\MTG(\k)$ devient $\pmt$, 
en vertu de l'égalité  tautologique
de $\Phi(\FC_1\FC_2)$ et $\FC_1(\Phi(\FC_2))$.

Le monoïde $\MTG$ est en fait un groupe. Pour cela, 
il suffit de vérifier la pleine fidélité
de tout $\FC\in\MTG(\k)$, qu'on peut se contenter de tester 
sur $\endo{\Pi(\k)}{O_0}$. Or $\kappa_G(\FC)$ est inversible, 
car pro-unipotent. 

L'ensemble des $\FC$ tels que le terme constant de 
$\Phi(\FC)$ soit 1 est clairement un sous-groupe de $\MTG(\k)$. Dans ce cas, 
l'action sur $\homp{O_0}{O_1}$ est elle-même pro-unipotente. Elle est 
fidèle car $\FC$ se déduit de $\FC(1_{1,0})$. Elle 
s'obtient par extension continue des scalaires de $\QM$ à $\k$. 

\begin{prop}
L'ensemble $\MT(\k)$ des séries de terme constant 1, muni de la loi 
$\pmt$ définit un groupe pro-unipotent agissant par $\pmt$ sur 
$\sernc{\k}{\XB}$.  
\end{prop}	

\sss{Remarques}\label{sss_remMT}
Par définition, l'action de $G\in\MT(\k)$ commute à la multiplication par
gauche par $x_1$ et à droite par $x_0$ et passe donc aux quotients par $\pi_Y$
et $\pi_\cv$. 

Si $S\in\MT(\k)$ 
ne dépend que de $x_1$, l'automorphisme $\kappa_S$ est l'identité,
car $t_\s(S)$, ne dépendant que de $x_\s$, commute à $x_\s$. L'action de 
$S$ est alors la multiplication à gauche par $S$ et on a 
$G\pmt S=G\pmt (S\pmt 1)=SG=S\pmt G$. Autrement dit, $S$ est central. 

\sss{}Soit $\MTD(\k)$ le sous-monoïde de $\MTG(\k)$ agissant sur les flèches
de $\Pi(\k)$ par morphismes de cogèbres. L'élément $\Phi(\FC)$ paramétrant 
$\FC\in\MTD(\k)$ est \glk{}, car $1_{1,0}$ l'est. Réciproquement, si $G$ 
est un élément \glk{} de $\sernc{\k}{\XB}$, le foncteur 
$\Phi^{-1}(G)$ appartient à $\MTD(\k)$, car son action sur 
$\homp{O_\s}{O_\t}$
est, comme le cas particulier (\ref{eq_defMT}), la composition d'une 
multiplication par un élément \glk{} et d'opérateurs du type $\kappa_G$, qui
sont tous des morphismes de cogèbres, comme on le voit par test sur les 
éléments de $\XB$. En résumé:
\begin{prop}\label{prop_expliestable}
   $(\exp(\liel{}{\XB}),\pmt)$ est un sous-groupe pro-unipotent 
de $\MT$ ; il agit par morphisme de cogèbres sur $(\sernc{\k}{\XB},\D)$. 
\end{prop}

\vdeux{}{
On peut également considérer les sous-groupes de $\MTG$ qui respectent 
d'autres symétries, par exemple, pour $G=\triv$,
l'action de $\SG_3$ sur $\Pi(\k)$ provenant de l'action sur 
$\Puntrois$, ou, pour $\G=\mub_n$, l'action de $\ZM/2$ 
sur $\Pi(\k)$ qui provient de l'inversion sur $\Punrac{n}$.} 

\sss{}\label{sss_MTfoncG}
Explicitons les dépendances en $\G$ par les notations $\Pi_\G$ et 
$\MT(\G)$. Un morphisme de groupes
$\phi:\G\to\G'$ donne un foncteur $\Pi(\G)\to\Pi(\G')$, agissant sur 
les objets par $O_\s\mapsto O_{\phi{\s}}$ et sur les flèches par 
$\phi^*$. Il commute aux $C_\s$. Ceci montre que $\phi^*$ est 
un morphisme de groupes pro-unipotents $\MT(\G)\to\MT(\G')$, respectant
les actions sur les $\sernc{\k}{\XB}$. 

Si $\phi$ est injectif, on définit un foncteur de la sous-catégorie
pleine de $\Pi(\G')$ d'objets les $O_{\phi(\s)}$ dans $\Pi(\G)$ 
par $O_{\phi(\s)}\mapsto O_\s$ et $\phi_*$ sur les flèches et commutant 
aux $C_\s$ appropriés. On en déduit 
que $\phi_*$ est un morphisme $\MT(\G')\to\MT(\G)$,
respectant les actions. 
\sss{Structure infinitésimale} 
L'algèbre de Lie $\mt$ du groupe $\MT$ est formée des séries de 
terme constant nul. On notera dans la suite $\Exp$ l'application exponentielle
$\mt\to\MT$, et $\ih{\pcdot}{\pcdot}$ le crochet de Lie\footnote{Cette notation
suit celle de Drinfel'd \cite[p. 851]{DrinQTQH}.} de $\mt$, pour les 
distinguer de l'exponentielle et du crochet usuels de $\sernc{\k}{\XB}$. 
Le crochet de $\mt$ sera appelé {\em crochet d'Ihara}. 
On peut expliciter cette structure grâce à la représentation dans 
$\sernc{\k}{\XB}$: 

Soit $\eps$ une variable formelle telle que $\eps^2=0$. 
Pour $\psi\in\mt(\k)$, notons $s_\psi$ l'opérateur qui associe à 
$v\in\sernc{\k}{\XB}$ le terme en $\eps$ de $(1+\eps\psi)\pmt v$. 
La machinerie standard sur les groupes pro-unipotents donne:
\begin{prop}
   Pour tous $\psi\in\mt(\k)$ et $H\in\sernc{\k}{\XB}$, 
   on a 
\begin{equation}
\Exp(\psi)\pmt H=\exp(s_\psi)(H) 
\end{equation}

Pour tous $\psi_1,\psi_2\in\mt(\k)$, on a 
\begin{equation}
s_{\ih{\psi_1}{\psi_2}}=[s_{\psi_1}, s_{\psi_2}]	
\end{equation}
\end{prop} 
Tenant compte de $s_{\psi}(1)=\psi$ et $G\pmt 1 = G$, on obtient en particulier:
\begin{cor} 
Pour 
tous $\psi_1,\psi_2\in \mt(\k)$, on a 
\begin{eqnarray}
\Exp(\psi)&=&\exp(s_\psi)(1)  \\
\ih{\psi_1}{\psi_2}&=&s_{\psi_1}(\psi_2)-s_{\psi_2}(\psi_1)\label{eq_ihs}
\end{eqnarray}  
\end{cor} 

\sss{}
À partir des définitions, il est facile d'expliciter les opérateurs $s_\psi$:
\begin{prop}
Pour tous $v\in\sernc{\k}{\XB}$ et $\psi\in\mt(\k)$, on a 
\begin{equation}
s_\psi(v) = \psi v + d_\psi(v),
\end{equation} 
où $d_\psi$ est la dérivation continue de $\sernc{\k}{\XB}$ caractérisée par: 
\begin{equation}
d_\psi(x_\s) = \left\{\begin{array}{ccc}
	[t_\s(\psi), x_\s]&\text{si}&\s\neq 0 \\
	0 &\text{si}&\s=0\end{array}\right.  
\end{equation} 
\end{prop}
On voit immédiatement que l'opération 
$(\psi_1,\psi_2)\mapsto d_{\psi_1}(\psi_2)$ est homogène pour le poids et
la longueur (\cf{} \ref{def_Xb}), 
ainsi donc que $(\psi_1,\psi_2)\mapsto s_{\psi_1}(\psi_2)$. 
L'expression (\ref{eq_ihs}) du crochet d'Ihara met en évidence 
la structure d'algèbre pré-Lie  de $\mt$ (voir \cite{Gerstcr}).

\Subsection{Énoncé détaillé}
\begin{theo}\label{theo_moi}
Fixons un plongement $\iota$ d'un groupe commutatif fini dans $\CM^*$.
\begin{itemize} 
   \item $\DMRP_0$ et $\DMRD_0$ sont des sous-schémas en groupes de $\MT$.  
   \item Pour tout $\QM$-anneau $\k$ et $\l\in\k$, les groupes 
	$\DMRP_0(\k)$ et $\DMRD_0(\k)$ agissent librement et transitivement
par multiplication $\pmt$ à gauche respectivement sur $\DMRP_\l(\k)$ et 
$\DMRD_\l(\k)$. 
   \item $\DMRP_\l(\k)$ et $\DMRD_\l(\k)$ sont non-vides. 
\end{itemize} 
\end{theo}
Autrement dit, $\DMRP\to\AM^1$ et $\DMRD\to\AM^1$ sont des torseurs triviaux,
respectivement sous $\DMRP_0$ et $\DMRD_0$.

\subsection{Espaces tangents à l'origine}\label{subsec_dm}
La première étape de la preuve du théorème \ref{theo_moi} concerne
les espaces tangents au voisinage de 1 de $\DMR$ et $\PDMRD$, \ie{} 
les points de la forme $1+\eps\psi$ à coefficients 
l'anneau des nombres duaux $\k[\eps]$. 

On les décrit en \ref{sss_defdm} et on étudie en 
\ref{sss_virelng1} les contraintes portant sur les termes de longueur 1
de leurs éléments. La propriété obtenue est l'origine combinatoire 
de l'intervention de la flèche $\DMRP\to\AM^1$. 
\sss{Définitions} \label{sss_defdm}
Pour tout $\QM$-anneau $\k$, et $\G$ quelconque, soit
$\dmr(\k)$ l'ensemble des séries $\psi$ de $\sernc{\k}{\XB}$ qui 
vérifient les équations: 
\begin{eqnarray}
(\psi|x_0)&=&(\psi|x_1)=0 \\
\D\psi=1\ot_\k\psi + \psi\ot_\k 1 &\text{et}& \Det(\psi_\et)=1\ot_\k\psi_\et +
\psi_\et\ot_\k 1,\label{eq_psiprim}\\
\qtext{o{\`u} l'on pose} \psi_\et\ \ass\ \qs\pi_Y(\psi) + \psi_\corr &\text{et}&
\psi_\corr\ \ass\ \sum_{n\geq 2}
\frac{(-1)^{n-1}}{n} (\psi|y_{n,1})y_{1,1}^n\nonumber
\end{eqnarray}

Soit $\dmrd(\k)$ l'ensemble de celles 
qui vérifient en outre, pour tout diviseur $d$ de l'ordre de $\G$, 
les relations de distribution 
\begin{equation}\label{dist_dmrd}
p^d_*(\psi) = i_d^*(\psi)+
\sum_{\s^{n/d}=1}(\psi|x_\s)x_1
\end{equation}

On laisse le lecteur se convaincre du fait que $\dmr$ et $\dmrd$ 
sont, au sens plus haut, les espaces tangents à $\DMR$ et $\PDMRD$ 
au voisinage de 1. 

Il est clair que la composante homogène de poids $p$ d'un élément de 
$\dmr(\k)$ est encore dans $\dmr(\k)$. On a donc là un schéma vectoriel
(\cf{} \ref{subsec_conv}), associé au $\QM$-espace vectoriel gradué, qu'on
notera encore $\dmr$, des éléments de $\assl{\QM}{X}$ qui sont dans 
$\dmr(\QM)$. De même pour $\dmrd$.

\sss{Termes de longueur 1}\label{sss_virelng1} On exhibe une contrainte portant sur les termes de longueur 1 des éléments de $\dmr$. On verra plus loin que 
c'est la seule. 
\begin{prop} \label{prop_virelng1}
Soient $\psi\in\dmr, \nu\in\G$ et $n$ un entier $\geq2$. On a les égalités: 
\begin{eqnarray}\label{eq_virelng1_a}
 \scal{\psi_\et}{y_{n,\nu}} + (-1)^n \scal{\psi_\et}{y_{n,\nu^{-1}}} &=& 0
\qtext{si $n\geq3$;}   \\ 
 \scal{\psi_\et}{y_{n,\nu}} + \scal{\psi_\et}{y_{n,\nu^{-1}}} + 
\scal{\psi_\et}{y_{2,1}} &=& 0\qtext{si $n=2$ et $\nu\neq1$.} 
\label{eq_virelng1_b} 
\end{eqnarray}
\end{prop}
\begin{lemme}\label{lemme_virelng1}
  Pour tout {\'e}l{\'e}ment $\psi$ de $\liel{}{\XB}$, tout
  $n\in\NM^*$ et tout $\nu\in\G$, on a   
$$ \scal{\qs\pi_Y(\psi)}{ y_{1,1}y_{n-1,\nu} + (-1)^n
  y_{1,\nu}y_{n-1,\nu^{-1}}} = 0$$
\end{lemme}
\begin{preuve} 
De mani{\`e}re {\'e}quivalente, on doit prouver que l'on a:
$$ \scal{\psi}{x_1x_0^{n-2}x_\nu + (-1)^n x_\nu x_0^{n-2}x_1} = 0$$
Comme $\psi\in\liel{}{X}$, l'antipode $S_X$ le transforme en son 
opposé. Le coefficient de $x_1x_0^{n-2}x_\nu$ est donc l'opposé de celui de 
$S_X(x_1x_0^{n-2}x_\nu)=(-1)^n x_\nu x_0^{n-2}x_1$
\end{preuve}
\begin{lemme}\label{lemme_virelng1_2} 
Soient $n\geq 2$ un entier, $\s,\tau\in\G$. 
 Si $n\geq 3$ ou $(\s,\t)\neq(1,1)$, on a, pour tout $\psi\in\liel{}{\XB}$:
\begin{equation}
 \scal{\psi_\et}{y_{n-1,\tau}y_{1,\s\tau^{-1}}} +
  \sum_{p=1}^{n-1} \scal{\psi_\et}{y_{p,\s}y_{n-p,\tau\s^{-1}}} = 0
\end{equation}
\end{lemme}
\begin{preuve}
On peut supposer $\psi$ homogène de poids $n$. 
Le terme correctif $\psi_\corr$ est dans ce cas un multiple 
scalaire de $y_{1,1}^n$. 
L'hypothèse sur $n,\s$ et $\tau$ assure donc que les coefficients de
$y_{1,\s}y_{1,\tau}$ dans $\qs\pi_Y(\psi)$ et $\psi_\et$ sont égaux.
En utilisant $\ps$, on est donc ramené à prouver:
\begin{equation}
   \scal{\psi}{x_0^{n-2} x_\tau x_\s} +
  \sum_{p=1}^{n-1} \scal{\psi}{x_0^{p-1} x_\s x_0^{n-p-1} x_\tau} = 0,
\end{equation}
Il s'agit en fait d'un calcul avec le produit $\sh$,
 dual de $\D$, dont on a jusqu'à présent évité l'emploi (voir \cite{Reut}). 
Il est défini par 
$\scal{a}{b\sh c}=\scal{\D(a)}{b\ot c}$ pour tous $a,b,c\in\assl{\QM}{\XB}$.  
Le développement de $x_\s \sh x_0^{n-2}x_\tau$ est 
la somme de toutes les insertions possibles de $x_\s$ dans $x_0^{n-2}x_\tau$:
$$x_\s \sh x_0^{n-2}x_\tau = 
x_0^{n-2} x_\tau x_\s +
  \sum_{p=1}^{n-1} x_0^{p-1} x_\s x_0^{n-p-1} x_\tau
$$
Comme $\psi$ est primitif, on a 
$\scal{\psi}{x_\s\sh x_0^{n-1}x_\tau}=\scal{\D\psi}{x_\s\ot x_0^{n-1}x_\tau}=0$
\end{preuve}
\sss{D{\'e}monstration de la proposition \ref{prop_virelng1}}
Le lemme \ref{lemme_virelng1_2}, appliqué, pour $n\geq3$ ou $\nu\neq1$, à 
$(\s,\tau)=(\nu, \nu)$, puis $(1,\nu)$ donne:
\begin{eqnarray}
\label{eq_virelng1_0} 
\scal{\psi_\et}{y_{n-1,\nu}y_{1,1}} + 
  \sum_{p=1}^{n-1}\scal{\psi_\et}{y_{p,\nu}y_{n-p,1}} &=& 0\qtext{et}\\
\scal{\psi_\et}{y_{n-1,\nu}y_{1,\nu^{-1}}} + 
  \sum_{p=1}^{n-1} \scal{\psi_\et}{y_{p,1}y_{n-p,\nu}} &=& 0,
	\qtext{\ie}\nonumber \\ 
\label{eq_virelng1_1} 
\scal{\psi_\et}{y_{n-1,\nu}y_{1,\nu^{-1}}} + 
   \sum_{p=1}^{n-1} \scal{\psi_\et}{y_{n-p,1}y_{p,\nu}} &=& 0
\end{eqnarray}
Or, pour tout $p$ compris entre 1 et $n-1$, la primitivité de $\psi_\et$
pour $\Det$ donne (\cf{} \ref{prop_etpart}): 
$$ \scal{\psi_\et}{y_{n,\nu} + y_{p,\nu}y_{n-p,1} + y_{n-p,1}y_{p\nu}} =
0 $$ On tire donc en sommant les {\'e}galit{\'e}s (\ref{eq_virelng1_0}) et
(\ref{eq_virelng1_1}): 
\begin{eqnarray}
  (\psi_\et | y_{n-1,\nu}y_{1,1} + y_{n-1,\nu}y_{1,\nu^{-1}}) -
  (n-1)(\psi_\et | y_{n,\nu}) &=& 0,\qtext{puis}\\
  -(\psi_\et | y_{1,1}y_{n-1,\nu} + y_{1,\nu^{-1}}y_{n-1,\nu}) -
  n(\psi_\et | y_{n,\nu}) - (\psi_\et | y_{n,1}) &=& 0 
\label{eq_virelng1_2}
\end{eqnarray}
L'{\'e}quation (\ref{eq_virelng1_2}) 
somm{\'e}e avec sa variante pour $\nu^{-1}$ pond{\'e}r{\'e}e par $(-1)^n$ 
donne:
\begin{multline}
n(\psi_\et | y_{n,\nu} + (-1)^n
y_{n,\nu^{-1}}) + (1 + (-1)^n) (\psi_\et | y_{n,1}) = \\
\left(\aug\psi_\et \left| y_{1,1}y_{n-1,\nu} + (-1)^n
y_{1,\nu}y_{n-1,\nu^{-1}}\right.\right) 
+ \left(\aug\psi_\et \left| y_{1,\nu^{-1}}y_{n-1,\nu} + (-1)^n
  y_{1,1}y_{n-1,\nu^{-1}}\right.\right)  
\end{multline}
D'apr{\`e}s le lemme (appliqu{\'e} pour $\nu$ et $\nu^{-1}$), on a donc
\begin{equation}\label{eq_virelng1_3}
 n(\psi_\et | y_{n,\nu} + (-1)^n
y_{n,\nu^{-1}}) + (1 + (-1)^n) (\psi_\et | y_{n,1}) = 0
\end{equation}

Pour $n$ impair et supérieur à $3$ ou $n=2, \nu\neq 1$, 
l'équation (\ref{eq_virelng1_3}) donne directement le résultat voulu.
Pour $n$ pair, supérieur à $3$ et $\nu=1$, l'équation 
(\ref{eq_virelng1_3}) devient $(\psi_\et|y_{n,1}) = 0$, ce qui 
reporté dans (\ref{eq_virelng1_3}) donne à nouveau (\ref{eq_virelng1_a}).
\hfill\qed
\sss{}
On a une réduction supplémentaire:
\begin{prop}\label{prop_y2} 
Soit $\psi\in\dmr$. Si $\G$ est de cardinal au moins 3, le coefficient de 
$y_{2,1}$ dans $\psi$ est nul.
\end{prop}
\begin{preuve}
\vdeux{Pour tout $\psi\in\assl{\QM}{\XB}$, }
{Pour tout $\psi\in\assl{\QM}{\XB}$, 
primitif pour $\D$, le coefficient de $x_\s\ot x_{\s\tau}$ dans 
$\D\psi$ étant nul, } on a  
 $\scal{\psi}{x_\s x_{\s\tau} + x_{\s\tau} x_\s} = 0$. 
   Si $(\s,\tau)\neq(1,1)$, on en déduit la nullité de  
$\scal{\ps\psi_\et}{x_\s x_{\s\tau} + x_{\s\tau} x_\s}$, ce qui donne:
$$
 (\psi_\et | y_{1,\s}y_{1,\tau}) = - (\psi_\et| y_{1,\s\tau}y_{1,\tau^-1}) 
= (\psi_\et | y_{2,\s} + y_{1,\tau^{-1}}y_{1,\s\tau}),$$
la deuxième égalité provenant de la proposition \ref{prop_etpart}.
L'application $(\s,\t)\mapsto(\t^{-1},\s\t)$ est cyclique d'ordre 6.
En itérant six fois l'égalité ci-dessus, on obtient facilement 
$$ \scal{\psi_\et}
{y_{2,\s} + y_{2,\tau^{-1}} + y_{2,\s^{-1}\tau^{-1}}
+ y_{2,\s^{-1}} + y_{2,\tau} + y_{2,\s\tau}} =0 $$
Grâce à la proposition \ref{prop_virelng1}, ceci entraîne la nullité de 
$\scal{\psi_\et}{y_{2,1}}$, à condition qu'il existe  
$\s$ et $\tau$ dans $\Gamma$ tels que $\s$, $\tau$ et $\s\tau$ soient tous
diff{\'e}rents de 1.

Un groupe $\G$ non trivial dans lequel un tel choix est impossible est 
de cardinal 2: 
fixons en effet $\s\in\G\setminus\{1\}$;
tout élément $\tau$ de $\G\setminus\{1\}$ est l'inverse de ce $\s$. 
\end{preuve}

\sss{Définitions de $\dmr_0$ et $\dmrd_0$}\label{sss_defdm0}
L'équation (\ref{eq_virelng1_a}) interviendra à plusieurs reprises comme 
condition technique dans la suite. C'est cette condition qui est à 
l'origine du fait qu'on obtienne une action fibre à fibre de $\DMRP_0$ sur
$\DMRP$, et non une structure de groupe sur $\DMRP$. 

On notera dans la suite $\dmr_0$ et $\dmrd_0$
l'ensemble des éléments $\psi$ de $\dmr$ (resp. $\dmrd$) 
qui vérifient (\ref{eq_virelng1_a}) pour tout $(n,\nu)\in\NM^*\times\G$. 
Par les propositions \ref{prop_virelng1} et \ref{prop_y2}, il suffit d'imposer
(\ref{eq_virelng1_a}) pour $n=1$ si $|\G|\geq3$ ; pour $(n,\nu)=(2,1)$ si 
$|\G|\leq2$. 

Pour tout plongement $\iota$ de $\G$ dans $\CM^*$, on voit que 
$\dmr_0$ est l'espace tangent à $\DMRP_0$  
au voisinage de 1. Il est indépendant de $\iota$, 
car les relations de poids 1 s'écrivent  
sans coefficients dans le cas dégénéré $\l=0$. 
De même $\dmrd_0$ est l'espace tangent à $\DMRD_0$ au voisinage de 1.

Comme cas particulier de (\ref{eq_virelng1_a}), on voit qu'un élément $\psi$ 
de $\dmr_0$ n'a pas de terme en $y_{n,1}$ 
si $n$ est pair. Son terme correctif $\psi_\corr$ est donc une série
impaire en $y_1$.

\section{Action tangente}\label{sec_tangente}
 On établit dans cette section les résultats ci-dessous, qui forment la 
partie la plus difficile de la démonstration du théorème. 
\begin{prop}\label{prop_DMRstable}
   Soient $\G$ un groupe commutatif fini et $\k$ un $\QM$-anneau.  
\begin{itemize}
	\item Les espaces tangents $\dmr_0$ et $\dmrd_0$ sont des 
sous-algèbres de Lie de $\mt$.
	\item Les ensembles 
$\DMR(\k)$ et $\PDMRD(\k)$ sont stables par $\exp(s_\psi)$, pour tout 
$\psi$ appartenant respectivement à $\dmr_0(\k)$ et $\dmrd_0(\k)$. 
	\item L'action de $\exp(s_\psi)$ sur $\DMR(\k)$ commute à
la flèche $\DMR(\k)\to\k$ définie par tout plongement de $\G$ dans $\CM^*$.  
\end{itemize}
\end{prop}
Grâce à la formule de Campbell-Hausdorff, le premier point implique que 
$\Exp(\dmr_0)$ et $\Exp(\dmrd_0)$ sont des sous-schémas en groupes 
de $\MT$. Le deuxième indique alors qu'il agit par multiplication $\pmt$ 
à gauche sur $\DMR$ (resp. $\PDMRD$). 

De \ref{subsec_infY} à \ref{subsec_DMRst}, 
on ne considère comme coproduit que
$\Det$. Par \og{}primitif\fg{}, il faut donc entendre par exemple 
\og{}primitif pour $\Det$\fg{}. 

\Subsection{Les opérateurs infinitésimaux de $\MT$, vus sur $\qyg$}\label{subsec_infY}
\sss{Définitions}
Pour tous $\psi,v\in\qxg$, comme $d_\psi$ est une dérivation qui annule 
$x_0$, on a 
$$s_\psi(vx_0)=\psi vx_0+d_\psi(vx_0)=\psi vx_0+d_\psi(v)x_0$$ 
Il s'ensuit que le noyau $\qxg x_0$ de $\pi_Y$ est 
stable par $s_\psi$. Soit alors $s_\psi^Y$, l'endomorphisme du $\k$-module 
$\qyg$ tel que le diagramme ci-dessous commute: 
$$\xymatrix{
\qxg \ar[r]^{s_\psi}\ar[d]^{\qs\pi_Y} & \qxg\ar[d]^{\qs\pi_Y} \\
\qyg \ar[r]^{s_\psi^Y} & \qyg}
$$

On considèrera également l'endomorphisme $D_\psi^Y$ de $\qyg$ donné par
$$D_\psi^Y(v)\ \ass\ s_\psi^Y(v) - v\qs\pi_Y(\psi)$$ 
\sss{} 
Pour calculer explicitement $s_\psi^Y(v)$, on doit choisir
$w\in\qxg$ tel que $\pi_Y(w)=v$. Comme l'inclusion de $\qyg$ dans $\qxg$ 
est une section de $\pi_Y$, on peut prendre $w=v$. 

Par exemple, on a $\pi_Y(x_1)=y_{1,1}$,
$\qs(y_{1,1})=y_{1,1}$ et $s_\psi(x_1)=\psi x_1 + [x_1,\psi]=x_1\psi$, d'où
on tire $\s_\psi^Y(y_{1,1})=\qs\pi_Y(x_1\psi)=y_{1,1}\qs\pi_Y(\psi)$. 
En d'autres termes:
\begin{prop}\label{prop_DpsiYy1}
Pour tout $\psi\in\qxg$, on a $D_\psi^Y(y_{1,1})=0$.
\end{prop}
On constate également la propriété suivante:
\begin{prop}\label{prop_Dpsider}
   Pour tout $\psi\in\qxg$, l'opérateur $D_\psi^Y$ est une dérivation de 
$\qyg$.
\end{prop}
\begin{preuve}
   Par linéarité et une récurrence évidente, il suffit de démontrer:
$$ D_\psi^Y(y_{n,\nu}w) = D_\psi^Y(y_{n,\nu})w + y_{n,\nu}D_\psi^Y(w),$$
pour tout $(n,\nu)\in\NM^*\times\G$ et tout $w\in\assl{\QM}{\YB}$, homogène
de couleur totale $\ov{w}$. 

Remarquons d'abord que, pour $v\in\assl{\QM}{\YB}$, homogène
de couleur totale
$\ov{v}$ on a $\ps(vw) = \ps(v)t_{\ov{v}}(w)$, ce qui se vérifie facilement dans le cas où $v$ et $w$ sont des mots de $\YB$ et s'écrit encore:
\fnote{Rappel : $t_\nu$ est la substitution $x_0\mapsto x_0, x_\s\mapsto 
x_{\s\nu}$.}
\begin{equation}\label{eq_deriv_r}
vw = \qs(\ps(v)t_{\ov{v}}\ps(w))
\end{equation}
D'après la définition de $\ps$, on voit d'autre part que $\qs(v)$ est 
homogène de couleur totale $\ov{v}$ si et seulement si tous les mots de
$\XB$ intervenant dans $v$ se terminent par $x_v$. 

En appliquant la définition de $D_\psi^Y$, on a
\begin{gather}
D_\psi^Y(y_{n,\nu}w) = \qs\pi_Y\left[\aug 
d_\psi(y_{n,\nu}t_\nu\ps(w)) + \psi y_{n,\nu} t_\nu \ps(w)\right] 
- y_{n,\nu}w\qs\pi_Y(\psi) \nonumber\\
 = \qs\pi_Y\left[\aug
d_\psi(y_{n,\nu})t_\nu\ps(w)) + \psi y_{n,\nu} t_\nu \ps(w) 
-y_{n,\nu}t_\nu\ps(w)t_{\nu\ov w}\psi \right]\nonumber\\
=  \qs\pi_Y\left[\aug 
d_\psi(y_{n,\nu})t_\nu\ps(w) + y_{n,\nu}t_\nu d_\psi(\ps(w))
+ \psi y_{n,\nu} t_\nu \ps(w) -y_{n,\nu}t_\nu\ps(w)t_{\nu\ov w}\psi \right],
\label{eq_deriv1}
\end{gather}
cette dernière ligne s'obtenant en écrivant que $d_\psi$ est une dérivation
et qu'elle commute, par définition, à $t_\nu$. 
On obtient de même 
\begin{equation}\label{eq_deriv2}
D_\psi^Y(y_{n,\nu}) = \qs\pi_Y\left[\aug
d_\psi(y_{n,\nu}) + \psi y_{n,\nu} -y_{n,\nu}t_\nu(\psi)\right]
\end{equation}
En revenant à la définition de $d_\psi$, on en déduit immédiatement que 
$D_\psi^Y(y_{n,\nu})$ est égal à 
$\qs\pi_Y[-x_0^{n-1}t_\nu(\psi)x_\nu+\psi y_{n,\nu}]$ et est donc 
homogène de couleur totale $\nu$. On peut donc multiplier
(\ref{eq_deriv2}) à droite par $\ps(w)$ en utilisant (\ref{eq_deriv_r}).
Il vient:
\begin{equation}
D_\psi^Y(y_{n,\nu}) w = \qs\pi_Y\left[\aug
d_\psi(y_{n,\nu})t_\nu\ps(w) + \psi y_{n,\nu}t_\nu\ps(w) 
-y_{n,\nu}t_\nu(\psi)\t_\nu(\ps(w))
\right]\label{eq_deriv3}
\end{equation}
D'un autre côté, on a, plus facilement 
\begin{eqnarray}
\nonumber y_{n,\nu}D_\psi^Y(w) &=& y_{n,\nu}\qs\pi_Y\left[\aug
d_\psi(\ps(w)) + \psi\ps(w) -\ps(w) t_{\ov w}(\psi)\right]\\
&=& \qs\pi_Y\left[\aug
y_{n,\nu}t_\nu d_\psi(\ps(w)) + y_{n,\nu}t_\nu(\psi\ps(w)) 
- y_{n,\nu}t_\nu(\ps(w))t_{\nu\ov w}(\psi)\right]
\label{eq_deriv4}
\end{eqnarray}
La somme des seconds membres de (\ref{eq_deriv3}) et 
(\ref{eq_deriv4}) est bien celui de (\ref{eq_deriv1}). 
\end{preuve}

\sss{Remarques} Dans le cas $\G=\triv$ des polyzêtas, l'opérateur 
$D_\psi^Y$ est simplement la restriction à $\qyg$ de la dérivation continue 
$D_\psi$ définie par $D_\psi(x_0)=[x_0,\psi]$ et $D_\psi(x_1)=0$. 
Les propositions \ref{prop_DpsiYy1} et \ref{prop_Dpsider} 
sont alors évidentes. Dans le cas général, on peut toujours considérer 
$D_\psi^Y$ comme la restriction de $v\mapsto \qs d_\psi \ps(v) - \psi\qs(v)$ 
à $\qyg$, mais ceci n'est plus une dérivation de $\qxg$.
 
\subsection{Remontée}Si $\psi_1$ et $\psi_2$ sont deux éléments de 
$\qxg$ ayant la même image 
par $\pi_Y$, il n'est pas vrai en général que $s_{\psi_1}^Y$ et 
$s_{\psi_2}^Y$ 
soient égales. On décrit ici une section de $\pi_Y$ adaptée aux 
éléments de $\liel{}{X}$. 
\sss{}
Notons $\d_{x_0}$ la dérivée partielle par rapport à $x_0$ dans 
$\assl{\QM}{\XB}$, c'est-à-dire la dérivation qui envoie $x_0$ sur $1$ 
et les $(x_\s)_{\s\in\G}$ sur 0. Elle est homogène pour toutes les 
graduations de $\assl{\QM}{\XB}$ pour lesquelles les éléments de $\XB$ sont
homogènes, en particulier le poids et la longueur. 
On vérifie facilement
qu'elle commute à $\ps$ et que $\qyg$ est stable par $\d_{x_0}$.  

Soit $\sec$ l'application de $\qyg$ dans $\qxg$ définie par 
\begin{equation}
\sec(\psi)\ \ass\ \sum_{i\geq 0} \frac{(-1)^i}{i!}\d_{x_0}^i(\psi)x_0^i
\end{equation}

Pour $\psi\in\qyg$, on a clairement $\pi_Y\sec\psi=\psi$. Tout élément 
$\psi$ de $\qxg$ s'écrit de manière unique sous la forme 
$\sum_{i\geq 0}\psi_i x_0^i$ si l'on impose aux $\psi_i$ 
d'appartenir à $\qyg$, et $\psi_0$ n'est autre que $\pi_Y(\psi)$. 
La condition $\d_{x_0}(\psi)=0$ est donc équivalente à 
$$ \sum_{i\geq 0}\d_{x_0}(\psi_i) x_0^i + 
\sum_{i\geq 0} (i+1)\psi_{i+1} x_0^i $$ 
Les $\d_{x_0}(\psi_i)$ appartenant à $\qyg$, ceci est encore équivalent à 
la relation de récurrence $(i+1)\psi_{i+1}+\psi_i=0$, 
pour tout $i\in\NM$. On en déduit 
que les conditions $\d_{x_0}(\psi)=0$ et $\psi=\sec(\psi_0)$ sont équivalentes.
Autrement dit:
\begin{prop}\label{prop_sec}
   L'application $\sec:\qyg\to\ker\d_{x_0}\subset\qxg$ est l'inverse de 
la restriction de $\pi_Y$ à $\ker\d_{x_0}$. 
\end{prop}

De la définition de $\d_{x_0}$ découle par une récurrence immédiate que 
$\d_{x_0}$ annule tout élément de $\liel{}{X}$ homogène
de poids au moins 2. Pour $\psi\in\whliel{\k}{\XB}$,
on a donc simplement $\d_{x_0}(\psi)=(\psi|x_0)$. De plus, on a 
$\d_{x_0}(y_{1,1}^n)=0$, et donc $\sec(y_{1,1}^n)=y_{1,1}^n$. De tout 
cela s'ensuit:
\begin{prop}\label{cor_sec}
   Pour $\psi\in\dmr(\k)$, on a $\psi=\ps\sec(\psi_\et) - \psi_\corr$. 
\end{prop}

\sss{Propriétés de $\d_{x_0}$ sur $\assl{\QM}{Y}$.}
La restriction de $\d_{x_0}$ à $\qyg$ est la dérivation
qui envoie $y_{n,\nu}=x_0^{n-1}x_\nu$ sur $(n-1)y_{n-1,\nu}$.
Elle est donc homogène de degré $(-1,1)$ 
pour le poids coloré et de degré -1 pour la longueur. 
\begin{prop}
   Pour tous $(n,\nu)\in\NM^*\times\G$, on a 
\vdeux{$\d_{x_0} u_{n,\nu}=(n-1)y_{n,\nu}$}
{$\d_{x_0} u_{n,\nu}=(n-1)u_{n,\nu}$}.
\end{prop} 
\begin{preuve}
Dans l'algèbre $\sernc{\QM[t]}{\YB}$, considérons la dérivation 
$\QM$-linéaire $\d_t$ par rapport à $t$. Notons $\d=\d_{x_0}-t^2\d_t$. 
Comme $t^2$ est central, $\d$ est encore une dérivation. 

Soit $(\psi_n)_{n\geq 0}$ une suite
d'éléments homogènes de poids $n$ 
de $\assl{\QM}{\YB}$ et 
\vdeux{$\Psi=\sum_{n}\psi_n$}{$\Psi\ \ass\ \sum_{n}\psi_nt^n$}.
On vérifie facilement que les conditions
\begin{equation}\label{eq_dx0un}
\left[\text{Pour tout}\ n>0,\ \d_{x_0}(\psi_n)=(n-1)\psi_{n-1}\right]
\qtext{et}\d\Psi=0
\end{equation}
sont équivalentes. La première est vérifiée pour 
$\psi_n=\sum_{\nu\in\G} y_{n,\nu}$, et dans ce cas, la série 
$\Psi$ n'est autre que $\YC_1(t)$, avec les notations de \ref{sss_sommepart}. 

Soit $\UC(t)$ la série $\sum_{n,\nu} t^n u_{n,\nu}$. Par définition 
des $u_{n,\nu}$, on a $\UC(t)=\log(\YC_1(t))$. 
Comme  $\YC_1(t)$ vérifie (\ref{eq_dx0un}), 
en développant le logarithme et en utilisant que $\d$ est une dérivation,
on en déduit $\d\UC(t)=0$. 

Pour tout $n>0$, on a donc $\d_{x_0}(\psi_n)=(n-1)\psi_{n-1}$, avec 
$\psi_n=\sum_{\nu\in\G}u_{n,\nu}$. La proposition en découle, car $\d_{x_0}$
est homogène pour la couleur totale. 
\end{preuve}

Comme l'ensemble des
éléments primitifs de $\assl{\QM}{\YB}$ est l'ensemble des polynômes de Lie
en les $u_{n,\nu}$, cela implique: 
\begin{cor}\label{cor_derAprim}
   L'ensemble des primitifs de $\qyg$ est stable par $\d_{x_0}$. 
\end{cor} 
On aurait aussi pu vérifier directement par calcul sur les $y_{n,\nu}$ que 
$\d_{x_0}$ est une codérivation, redonnant ainsi ce corollaire, 
\subsection{Action par codérivation}
On établit ici par calcul direct le résultat ci-dessous, dont la proposition 
\ref{prop_DMRstable} sera une conséquence quasi-directe. 
\begin{prop} \label{prop_coder}
  Pour tout $\psi\in\dmr_0$, l'application 
$s_{\ps\sec(\psi_\et)}^Y$ est une codérivation. 
\end{prop}

\sss{Retournements} 
Notons $\ret_Y$ l'anti-automorphisme d'algèbres de $\qyg$ tel que 
$ret_Y(y_{n,\nu})=y_{n,\nu^{-1}}$, pour $(n,\nu)\in\NM^*\times\G$. 
Par test sur les $y_{n,\nu}$, c'est un morphisme de cogèbres qui fixe 1,
d'où en particulier:
\begin{prop}\label{prop_retYprim}
L'ensemble des éléments primitifs de $\qyg$ est stable 
par $\ret_Y$. 
\end{prop}

Cet opérateur apparaît dans l'expression des $D_\psi^Y(y_{n,\nu})$ lorsque 
$S_X(\psi)=-\psi$. 
\footnote{Rappelons que $S_X$, l'antipode de $(\assl{\QM}{\XB},\cdot\D)$, 
est l'anti-automorphisme d'algèbres 
 de $\assl{\QM}{\XB}$ qui envoie $x_\s$ sur $-x_\s$, pour tout $x_\s\in\XB$}
Cette condition est vérifiée en particulier par les éléments de 
$\liel{}{\XB}$. 
\begin{lemme}\label{lemme_Dpsiyn}
   Soit $\psi$ un élément homogène de poids $p$ de $\assl{\k}{\XB}$ tel que  
   $S_X(\psi)=-\psi$. Soient $(\psi_{i,\g})_{0\leq i\leq p, 
\g\in\G}$ les éléments de $\qyg$ de couleur totale $\g$ caractérisés par:
\begin{equation*}
   \psi = \sum_{i\geq0}\ps(\psi_{i,\g})x_0^i
\end{equation*}
   Avec ces notations, on a pour tout $(n,\nu)\in\NM^*\times\G$:
\begin{equation}\label{eq_Dpsilng1}
   D_\psi^Y(y_{n,\nu}) = \sum_{i\geq 0}\left(\aug\psi_iy_{n+i,\nu\gamma^{-1}}
+(-1)^py_{n+i,\nu\gamma}\ret_Y(\psi_{i\gamma})\right)
\end{equation}
\end{lemme}
\begin{preuve}
    On a successivement 
\begin{eqnarray*}
s_\psi(y_{n,\nu})&=&s_\psi(x_0^{n-1}x_\nu)=\psi x_0^{n-1}x_\nu + x_0^{n-1}x_\nu t_\nu(\psi) - x_0^{n-1}t_\nu(\psi)x_\nu, \qtext{d'où} \\
s_\psi^Y(y_{n,\nu})&=& \qs(\psi x_0^{n-1}x_\nu) + y_{n,\nu}\pi_Y\psi 
- \qs(x_0^{n-1}\psi x_\nu)  
\end{eqnarray*}
Le résultat voulu découle alors de la définition de $D_\psi^Y$ et du 
lemme immédiat ci-dessous, compte-tenu de $S_X(\psi)=-\psi$. 
\end{preuve}
\begin{lemme}\label{lemme_DpsiYn0}
  Pour tout {\'e}l{\'e}ment $\psi=\sum_{i\geq0,\g\in\G}\ps(\psi_{i,\g})x_0^i$ de $\qxg$, 
homog{\`e}ne de poids $p$, les  $\psi_{i,\g}$ {\'e}tant homogènes de couleur 
totale $\g$ dans $\qyg$ et tout $(n,\nu)\in\NM^*\times\G$, les formules 
suivante sont valables: 
\begin{eqnarray*} 
\qs(\psi x_0^{n-1}x_\nu)
   &=&\sum_{i\geq0, \g\in\G}\psi_{i,\g} y_{n+i,\nu\gamma^{-1}}\\
\qs(x_0^{n-1}t_\nu S_X(\psi)x_\nu)
  &=&(-1)^p\sum_{i\geq0}{y_{n+i,\nu\gamma}\ret_Y(\psi_{i,\g})}
\end{eqnarray*}
\end{lemme}
\sss{Abréviations}\label{sss_abrev}
Étant donné que $s^Y_{\ps\sec(\psi_\et)}$ dépend linéairement de $\psi$, pour 
démontrer la proposition \ref{prop_coder}, il suffit de traiter le cas où
$\psi$ est un élément homogène de poids $p$ de $\dmr_0$. On le supposera
fixé dans la suite et on adoptera les notations suivantes:
\begin{itemize}
   \item on abrège $D_{\ps\sec(\psi_\et)}^Y$ en $D$
   \item pour tout $(i,\g)\in\NM\times\G$, on note $\psi_{i,\gamma}$ la 
composante homogène de couleur totale $\g$ de $((-1)^i\d_{x_0}^i/i!(\psi_\et)$
et on pose 
$\chi_{i,\g}\ \ass\ (-1)^p \ret_Y(\psi_{i,\g})$. 
   \item on pose $z_{i,k}^{\g,\nu}\ \ass\ \psi_{i,\g}y_{k,\nu\g^{-1}} + y_{k, \nu\g}\chi_{i,\g}$, pour tout $(i,k,\g,\nu)\in\NM^2\times\G^2$.   
\end{itemize}
D'après la stabilité des primitifs par $\d_{x_0}$ et $\ret_Y$ 
(cor. \ref{cor_derAprim} et prop. \ref{prop_retYprim}) et l'homogénéité de 
$\Det$ pour la couleur totale, comme $\psi_\et$ est primitif, 
tous les $\psi_{i,\g}$ et $\chi_{i,\g}$ le sont. 

\sss{Calcul de $D$ sur les $y_{n,\nu}$}
D'après le corollaire \ref{cor_sec}, 
$\ps\sec(\psi_\et)$ est la somme du polynôme de Lie $\psi$ de 
$\qxg$ et de $\psi_\corr$, par définition multiple scalaire de 
$\scal{\psi}{y_p}x_1^p$, et donc nul par la proposition \ref{prop_virelng1} 
si $p$ est pair. 
Il s'ensuit que $\ps\sec(\psi_\et)$ est transformé en son opposé par $S_X$, 
ce qui permet d'appliquer le lemme \ref{lemme_Dpsiyn}. Avec les abréviations 
de \ref{sss_abrev}, cela s'écrit, pour $n>0$ et $\nu\in\G$:
\begin{equation}\label{eq_propcoder05}
D(y_{n,\nu})=\sum_{(i,\g)\in\NM\times\G} z^{\g,\nu}_{i,i+n}
\end{equation}
Ceci reste en fait vrai pour $n=0$, avec la convention (\ref{conv_y0}):
c'est le cas particulier $k=0$ du lemme ci-dessous, qui interviendra 
également dans le calcul final. 
\begin{lemme}\label{lemme_final}
Pour tout entier  $(k,\kappa)\in\NM\times\G$, on a 
\begin{equation}\label{eq_Dcoder2}
\sum_{i\geq k, \gamma\in\G} z^{\gamma,\kappa}_{i,i-k} = 0 
\end{equation}
\end{lemme}
\begin{preuve}
On a $D(y_{1,1})=0$ (prop. \ref{prop_DpsiYy1}). Avec la 
 formule (\ref{eq_propcoder05}) cela s'écrit
$$ \sum_{i\geq 0, \gamma\in\G} z^{\g,1}_{i,i+1} = 0 $$
D'après la proposition \ref{prop_deru}, en appliquant la dérivation 
$\partial_{u_{k+1,\kappa^{-1}}}$ à cette égalité, on obtient:
$$ \sum_{i\geq k, \gamma\in\G}
\left(\aug y_{i+1,\g^{-1}}(\psi_{i,\g}|y_{k+1,\kappa^{-1}}) + 
y_{i+1,\g}(\chi_{i,\g}|y_{k+1,\kappa^{-1}}) + z_{i,i-k}^{\gamma,\kappa}
\right) = 0
$$
Pour achever la démonstration du lemme, Il suffit donc de prouver l'égalité 
$$ \sum_{i\geq k, \gamma\in\G}
\left(\aug y_{i+1,\g^{-1}}(\psi_{i,\g}|y_{k+1,\kappa^{-1}}) + 
y_{i+1,\g}(\chi_{i,\g}|y_{k+1,\kappa^{-1}})\right)=0 $$ 
Les $\psi_{i,\g}$ et $\chi_{i,\g}$ étant homogènes de poids $p-i$ et 
de couleurs totales respectives $\g$ et $\g^{-1}$, les termes de cette somme 
autres que
$(i,\g)=(p-k-1,\kappa)$ et 
$(i,\g)=(p-k-1,\kappa^{-1})$ sont nuls ; celle-ci 
 vaut donc:  
\begin{eqnarray*}
&& y_{p-k,\kappa^{-1}}(\psi_{p-k-1,\kappa^{-1}}|y_{k,\kappa^{-1}})
+ y_{p-k,\kappa^{-1}}(\chi_{p-k-1,\kappa^{-1}}|y_{k,\kappa^{-1}})\\
&=& \frac{(p-1)!}{k!}y_{p-k,\kappa^{-1}}(\psi_\et | y_{p,\kappa^{-1}}
+ (-1)^p y_{p,\kappa}),
\end{eqnarray*}
et ceci est nul par définition de $\dmr_0$. 
\end{preuve}

\sss{Démonstration de la proposition \ref{prop_coder}}
La multiplication à droite par $\psi_\et$ est une codérivation car 
$\psi_\et$ est primitif. Il suffit donc de démontrer que 
$D=D_{\ps\sec(\psi_\et)}^Y$ est 
une codérivation. Comme $D$ est une dérivation, il suffit d'obtenir 
l'identité de codérivation sur les générateurs de $\qyg$, \ie{}
\begin{equation}\label{eq_Dcoder}
\forall(n,\nu)\in\NM^*\times\G,\ 
\Det D(y_{n,\nu}) = (\id\whot D + D\whot\id)\Det(y_{n,\nu}) 
\end{equation}

Abrégeons $a\ot b + b\ot a$ en $\sym(a\ot b)$.
La primitivité des $\psi_{i,\g}$ et $\chi_{i,\g}$ permet d'obtenir 
\begin{equation}\label{eq_Detz}
\Det(z^{\g,\nu}_{i,j}) = \underset{\kappa\l=\nu}{\sum_{k+l=i}} 
\sym(y_{k,\kappa}\ot z^{\g,\l}_{i,l})
\end{equation}
Le premier membre de (\ref{eq_Dcoder}) est donc, en utilisant la formule 
(\ref{eq_propcoder05}): 
\begin{equation}
  \Det D(y_{n,\nu}) = 
\underset{\g\in\G}{\sum_{i\geq 0}}\,
\underset{\kappa\in\G}{\sum_{k=0}^{i+n}}
\sym\left(\aug z_{i,k}^{\g,\kappa}\ot y_{n+i-k, \nu\kappa^{-1}}\right)
\end{equation}

Évaluons le second membre de (\ref{eq_Dcoder}):
$$\begin{array}{ccccc}
&&(D\ot\id+\id\ot D)\Det(y_{n,\nu}) &=& (D\ot\id +\id\ot D)
 \underset{\kappa\l=\nu}{\sum\limits_{k+l=n}} 
y_{k,\kappa}\ot y_{l,\lambda}\\
&=& \underset{\kappa\l=\nu}{\sum\limits_{k+l=n}} 
\sym\left(\aug D(y_{k,\kappa})\ot y_{l,\lambda}\right)
&=& \underset{\kappa\in\G}{\sum\limits_{k=0}^n}\, 
\underset{\gamma\in\G}{\sum\limits_{i\geq 0}}
\sym\left(\aug z_{i,i+k}^{\gamma,\kappa}\ot y_{n-k,\nu\kappa^{-1}}\right)
\\[5ex]
&=& \underset{\g\in\G}{\sum\limits_{i\geq 0}}\,
\underset{\kappa\in\G}{\sum\limits_{k=i}^{i+n}}
\sym\left(\aug z_{i,k}^{\g,\kappa}\ot y_{n+i-k, \nu\kappa^{-1}}\right)
\end{array}$$

La différence des deux membres de (\ref{eq_Dcoder}) est donc 
\begin{gather*}
\underset{\g\in\G}{\sum_{i\geq 0}}\,
\underset{\kappa\in\G}{\sum_{k=0}^{i-1}}
\sym\left(\aug z_{i,k}^{\g,\kappa}\ot y_{n+i-k, \nu\kappa^{-1}}\right)
= \underset{\g\in\G}{\sum_{i\geq 0}}\,
\underset{\kappa\in\G}{\sum_{k=1}^{i}}
\sym\left(\aug z_{i,i-k}^{\g,\kappa}\ot y_{n+k, \nu\kappa^{-1}}\right)\\
= \underset{\kappa\in\G}{\sum_{k\geq 1}}
\sym\left[\Aug
  \left(\aug \underset{\g\in\G}{\sum_{i\geq k}} z_{i,i-k}^{\g,\kappa}\right)
  \ot y_{n+k,\nu\kappa^{-1}}\right]
\end{gather*}
Or cette dernière expression est nulle d'après le lemme \ref{lemme_final}.\qed
\subsection{Preuve de la proposition \ref{prop_DMRstable}}
\label{subsec_DMRst}
Il nous reste essentiellement à vérifier que les termes correctifs se 
comportent bien.

\sss{}\label{calc_spsi1} Soient $\psi,v\in\qxg$. 
On a $d_\psi(x_1 v)=[x_1,\psi] v + x_1 d_\psi(v)$, d'où l'on tire 
\begin{equation}
   s_\psi(x_1v) = x_1 s_\psi(v)
\end{equation}
Soit $n$ un entier. La dérivation $d_{x_1^n}$ annule par définition $x_0$ et,
 pour tout $\s\in\G$, on a $d_{x_1}(x_\s)=[x_\s,x_\s^n]=0$. On en déduit donc
\begin{equation}
   s_{x_1^n}(v) = x_1^n v
\end{equation}
Il résulte de ces deux informations que $s_\psi$ et $s_{x_1^n}$ commutent. 
Tout ceci n'est que le cas particulier $\k=\QM[\eps]$ des remarques de 
\ref{sss_remMT}.

On a $\qs\pi_Y(x_1 v)=y_{1,1}\qs\pi_Y(v)$, d'où 
des formulations analogues pour les $s_\psi^Y$.
\sss{Termes de longueur 1}\label{ssslng1}
L'application bilinéaire $(\psi,\Phi)\mapsto s_\psi(\Phi)$ est homogène pour la longueur. Il en résulte si $\psi$ est sans 
terme constant que la partie de longueur 1 de $s_\psi^n(\Phi)$ est nulle pour
$n>1$. 
De plus, si le terme constant de $\Phi$ vaut 1, le terme en $y_{n,\nu}$ de 
$s_\psi(\Phi)$ est celui de $s_\psi(1)=\psi$. Ceci s'écrit également
\begin{equation}
 \scal{\exp(s_\psi)(\Phi)}{y_{n,\nu}} = \scal{\psi+\Phi}{y_{n,\nu}}
\end{equation}
De même, le terme de poids 1 de $\exp(s_\psi)\Phi$ est celui de $\psi+\Phi$. 

Cela s'applique notamment au cas où $\psi$ et $\Phi$ sont respectivement une 
série et une exponentielle de Lie. 

\sss{Crochet d'Ihara} 
Rappelons que l'on a $s_\psi(1)=\psi$, et donc $\qs\pi_Y(\psi)=s^Y_\psi(1)$,
pour tout $\psi\in\qxg$.

Soient $\psi_1$ et $\psi_2$ deux éléments de $\dmr_0$. L'algèbre de Lie 
libre $\liel{\k}{\XB}$ est stable pour le crochet d'Ihara (prop. \ref{prop_expliestable}). Celui-ci
étant homogène pour le poids et la longueur,
$\ih{\psi_1}{\psi_2}$ n'a
aucun terme de poids 1, et donc pas de terme en $x_0$ ni $x_1$, ni de 
longueur 1 et vérifie donc (\ref{eq_virelng1_a}) pour tout $n$. 

D'autre part, 
$s_{\ih{\psi_1}{\psi_2}}$ est le crochet de $s_{\psi_1}$ et $s_{\psi_2}$. 
Comme $\psi_i=\ps\sec\psi_{i,\et}+\psi_{i,\corr}$, pour $i=1,2$, on déduit de 
\ref{calc_spsi1} l'égalité 
\begin{eqnarray*}
 [s_{\psi_1},s_{\psi_2}] &=& [s_{\ps\sec\psi_{1,\et}}, s_{\ps\sec\psi_{2,\et}}],
\qtext{d'où}\\{}
[s^Y_{\psi_1},s^Y_{\psi_2}](1) &=& 
[s^Y_{\ps\sec\psi_{1,\et}}, s^Y_{\ps\sec\psi_{2,\et}}](1)
\end{eqnarray*}
Le membre de gauche de cette dernière égalité vaut 
$\qs\pi_Y(\ih{\psi_1}{\psi_2})$. Le membre de droite est l'image de 1 par une 
codérivation pour $\Det$: c'est donc un élément primitif pour $\Det$. 
On a donc prouvé que $\dmr_0$ est une sous-algèbre de Lie de $\mt$. 

\sss{Fin de la preuve pour $\DMR${}} Soit $\Phi\in\DMR(\k)$ et $\psi\in\dmr_0(\k)$. Il 
s'agit de prouver que $\Exp(\psi)\pmt\Phi=\exp(s_\psi)(\Phi)$ 
appartient à $\DMR(\k)$. 

Le coefficient de $x_0$ (resp. $x_1$) dans $\exp(s_\psi)(\Phi)$ est nul,
car c'est la somme de celui de $\psi$ et celui de $\Phi$. On sait déjà
(prop. \ref{prop_expliestable}) que $\exp(s_\psi)(\Phi)$ est 
\glk{}\vdeux{}{ pour $\D$}.

L'exponentielle d'une codérivation est un morphisme de cogèbres. 
L'élément 
$$G\ \ass\ \exp(s^Y_{\ps\sec(\psi_\et)})(\Phi_\et)$$ 
est donc \glk{} pour $\Det$. Comme $\ps\sec(\psi_\et)$ et $\psi$ (resp.
$\Phi_\et$ et $\qs\pi_Y(\Phi)$) sont égales à 
l'addition (resp. la multiplication à gauche) d'une série ne dépendant que 
de $x_1 (=y_{1,1})$ près, on déduit 
des remarques de \ref{calc_spsi1} que $G$ et 
$\qs\pi_Y\exp(s_\psi)(\Phi)$ sont égales, à multiplication à gauche près 
par une série ne dépendant que de $y_{1,1}$. Celle-ci est 
alors automatiquement $\left[\exp(s_\psi)(\Phi)\right]_\corr$, d'après
les résultats de \ref{sss_relregul}. 
%
%


\sss{Distribution}
Soient $\Phi_1,\Phi_2$ deux séries de termes constants 1 et 
satisfaisant, pour tout 
diviseur $d$ de l'ordre de $\G$ à (\ref{dist_DMRD}): 
on a $S_k\in\serc{\k}{x_1}$ tel que 
$p^d_*(\Phi_k)=S_k\cdot i_d^*(\Phi_k)$ pour $k\in\{1,2\}$.

Les applications $i_d^*$ et 
$p^d_*$ sont des morphismes de schémas en groupes de $\MT(\G)$ 
dans $\MT(\G^d)$. Comme les $S_k$ sont de plus centrales dans $\MT(\G^d)$, 
la relation de distribution vaut pour $\Phi_1\pmt\Phi_2$, à multiplication
à gauche près par $S_1S_2$, qu'on calcule facilement. Si les $\Phi_k$ sont 
des exponentielles de Lie, sans termes en $x_0$ ni $x_1$, 
on conclut plus rapidement par les arguments de \ref{sss_distrib}. 

Si $\psi$ est une série de Lie vérifiant (\ref{dist_dmrd}),  
$p^d_*(\psi)=S+i_d^*(\psi)$, où $S$ ne dépend que de $x_1$, il 
est clair pour les mêmes raisons que $\Exp\psi$ vérifie (\ref{dist_DMRD}), 
si les coefficients de $x_0$ et $x_1$ dans $\psi$ sont nuls. 
\sss{Poids un}
Soit $\iota$ un morphisme de groupes $\G\to\CM^*$. L'élément $\alpha_\iota$ 
de la définition \ref{def_DMl} est homogène de longueur 1. 
Par définition, si $\psi$ appartient à $\dmr_0(\k)$, on a 
$\scal{\psi}{\alpha_\iota}=0$. D'après \ref{ssslng1}, on a donc pour 
tout $\Phi\in\sernc{\k}{\XB}$:
$$\scal{\exp(s_\psi)(\Phi)}{\alpha_{\iota}}=\scal{\Phi}{\alpha_\iota}$$
Ceci donne la stabilité des relations de poids un et le point 
iii) de la proposition \ref{prop_DMRstable}.

\section{Transitivité}\label{sec_transitivite}
 \def\Mm{\AM^1}
\def\Y{\wh{C}}
\def\Gal{\text{\rm Gal}}
\Subsection{Approximations successives}\label{subsec_approx}
\sss{}
Il est classique d'interpréter les actions 
de $\Gm$ en termes $\ZM$-graduations. 
De même, 
une action du monoïde multiplicatif $\Mm$ sur un schéma affine 
$X=\spec A$ correspond à la donnée 
d'une $\NM$-graduation sur $A$: on munit $A\times\QM[t]$ de la graduation 
portée par le membre de droite, et on la ramène à $A$ par le morphisme 
d'algèbres \vdeux{$A\to A\times\QM[t]$}{$A\to A\ot\QM[t]$} 
définissant l'action.\fnote{Préciser le 
rôle de l'associativité de l'action}.

On supposera dans tout ce qui suit que 
les composantes homogènes de $A$ sont de dimension finie. 
La multiplication de $A$ fait alors 
du dual gradué de $A$ une cogèbre graduée 
$(C,\eps, \D)$ dont on notera $C_n$ la composante homogène de degré $n$. 
Étendons $\D$ et $\eps$ au complété $\Y(\k)$.
L'ensemble des $\k$-points de $X$ s'interprète alors 
comme l'ensemble des éléments \glks{} de 
$(\Y(\k), \eps, \D)$. Dans ce cadre, $\l\in\k=\Mm(\k)$ agit sur la composante 
homogène de degré $n$ de $\Y(\k)$ par multiplication par $\l^n$, et ceci
redonne l'action de $\Mm$ sur $X$. 

\sss{} Identifions le $(n+1)$\eme{} quotient $\tr{\Y}{n}(\k)$ de 
$\Y(\k)$ par la filtration associée à la graduation de $C$ 
à $\oplus_{i=0}^n C_i\ot\k$ et notons $\tr{\pi}{n}$ la 
projection de $\Y(\k)$ sur $\tr{\Y}{n}(\k)$ correspondante. 
 
Soit $\tr{X}{n}(\k)$ l'ensemble des éléments $\Phi$ de 
$\tr{\Y}{n}(\k)$ qui sont 
\glks{} modulo des termes de degré $n+1$, \ie{} qui vérifient
\begin{equation}\label{def_Xn}
\D\Phi = (\tr{\pi}{n}\whot\tr{\pi}{n})(\Phi\whot\Phi)\qtext{et}\eps(\Phi)=1
\end{equation}
Il est clair que $X$ est la limite projective des $\tr{X}{n}$. 

\sss{}
Par définition, un $\k$-point de 
l'espace tangent $T_\Phi X$ à $X$ au voisinage de $\Phi\in X(\QM)$ 
est un élément $\psi$ de $Y(\k)$ tel que $\Phi+\eps\psi$ soit un 
$\k[\eps]$-point de $X$, \ie{} \glk{} dans $Y(\k[\eps])$. Cela se traduit par
la condition 
\begin{equation}\label{eq_tangent}
   \D\psi = \Phi\ot\psi + \psi\ot\Phi
\end{equation}
Si $\Phi$ est stable par l'action de $\Mm$, il en est de même de l'espace 
tangent, qui est donc un schéma vectoriel associé à un espace vectoriel 
gradué. Plus concrètement, les composantes homogènes 
de degré $>0$ de $\Phi$, vu comme élément de $\Y(\k)$, sont nulles. 
Par homogénéité de $\D$, les composantes homogènes 
d'un élément $\psi$ de $T_\Phi(X)$ sont encore dans $T_\Phi(X)$, et 
$\psi$ en est la somme infinie. Si l'on prefère voir $\Phi$ comme un 
morphisme d'algèbres $A\to\QM$, cela fait de $T_\Phi X(\k)$ l'ensemble 
des $\Phi$-dérivations de $A$ dans $\k$.  

\sss{} \label{sss_hyp}
Supposons $X$ muni d'un morphisme $X\To{\alpha}\AM^1$, dont 
on notera $X_\l$ la fibre au-dessus d'un point $\l\in\k$, jouissant des
propriétés suivantes:
\begin{enumerate}
\item Homogénéité: le diagramme ci-dessous est commutatif. 
\begin{equation}
 \xymatrix{X\times \Mm \ar[rr]\ar[d]_{\alpha\times\id} && X\ar[d]^{\alpha} \\
 \Mm\times\Mm \ar[rr]_-{(\l,\mu)\mapsto \l\mu} && \Mm}
\end{equation}
\item Il existe un élément de $X_0(\QM)$ 
stable par l'action de $\Mm$ et un seul. On le note $1$. 
\item Soit $\xG$ l'espace tangent à $X$ au voisinage de $1$. La fibre spéciale
$\xG_0$ est une algèbre de Lie et on a une action  
de $\exp(\xG_0)$ sur $X$, qui est homogène, 
\ie{} qui respecte les actions de $\Mm$ et commute à $\alpha$. 
\item L'application $\exp(\xG_0)(\k)\to X_0(\k)$ donnée par l'action 
sur $1$ est injective pour tout $\k$.
\item Il existe un $\QM$-anneau $\KM$ et $\l\in\KM$ inversible tel que 
$X_\l(\KM)$ soit non-vide.
\end{enumerate}

\begin{prop} \label{prop_trans}
Sous ces hypothèses, l'action de $\exp(\xG_0)$ sur $X$ est
transitive et $X_\l(\k)$ est non-vide pour tout $\k$ et tout $\l\in\k$. 
On a un isomorphisme $X_0\simeq\exp(\xG_0)$.
\end{prop}
Les deux paragraphes suivants sont consacrés à la preuve de cette proposition.
\sss{} 
Le morphisme $X\to\AM^1$ donne un morphisme de cogèbres 
$\Y(\k)\to\serc{\k}{t}$. Par l'hypothèse d'homogénéité, il commute 
aux $\tr{\pi}{n}$. 

La formule (\ref{def_Xn}) implique en particulier que le 
terme de degré $0$ d'un élément de $\tr{X}{n}(\k)$ est \glk{}. Il est 
donc égal à 1, par l'hypothèse ii). 

\begin{prop}\label{prop_etape}
   la différence de deux éléments de $\tr{X}{n+1}_\l(\k)$ ayant même 
image par $\tr{\pi}{n}$ est un élément homogène de degré $n+1$ de $\xG_0(\k)$.

Si un élément de $\tr{X}{n}(\QM)$ admet un relevé dans $\tr{X}{n+1}(\k)$, 
il en admet un dans \vdeux{$\tr{X}{n}(\QM)$}{$\tr{X}{n+1}(\QM)$}.
\end{prop}
\begin{preuve}
Écrivons les composantes homogènes d'un élément $\Phi\in\tr{X}{n}(\k)$: 
\begin{equation}
\Phi = 1+\Phi_1+\cdots+\Phi_n
\end{equation}
La condition portant sur un élément $\Phi_{n+1}$ de $Y(\k)$, homogène de degré $n+1$, pour que $\Phi+\Phi_{n+1}$ appartienne à $\tr{X}{n+1}$ s'écrit:
$$
\D(\Phi_{n+1})-1\ot\Phi_{n+1}-\Phi_{n+1}\ot 1
=\underset{k,l>0}{\sum_{k+l=n+1}}\Phi_k\ot\Phi_l
$$
Sous cette forme, la première assertion est évidente, compte-tenu de la 
caractérisation (\ref{eq_tangent}) de $\xG(\k)$ et de la linéarité de 
$\Y(\k)\to\serc{\k}{t}$. 
 La seconde revient à dire qu'un 
système linéaire à coefficients rationnels ayant une 
solution dans $\k$ admet une solution dans $\QM$. 
\end{preuve}

\sss{Prise en compte de l'action} Un $\k$-point de $\exp(\xG_0)$ agit 
sur $\Y(\k)$ par morphisme de cogèbres. En exprimant ceci pour les 
nombres duaux $\k[\eps]$, 
on trouve une action, homogène, de $\xG_0$ sur $\Y$ par codérivations
qu'on notera encore $(\psi, v)\mapsto s_\psi(v)$. 
À nouveau, l'action de $\Exp(\psi)$ sur $\Y(\k)$ se fait par $\exp(s_\psi)$, 
pour $\psi\in\xG_0(\k)$. 

L'application $\psi\mapsto s_\psi(1)$ est un endomorphisme 
$\k$-linéaire de la partie homogène de degré $n$ de $\xG_0(\k)$. 
L'hypothèse iv) appliquée à
$\QM[\eps]$ montre qu'il est injectif pour $\k=\QM$, 
donc bijectif par finitude de la dimension. 
Par extension des scalaires, il est donc inversible pour tout $\k$. 

Enfin, on définit pour tout $n$ une action de $\exp(\xG_0)$ sur 
$\tr{\pi}{n}(\Y)$ par l'action sur $\Y$, suivie de $\tr{\pi}{n}$. Avec
l'homogénéité et le fait que $s_\psi$ est une codérivation, 
on voit que $\tr{X}{n}$ est stable par cette action. 
Si $\psi$ est homogène de degré $n$, 
l'action de $\exp(s_\psi)$ sur $\tr{X}{n}$ est simplement l'addition de 
$s_\psi(1)$. 

\begin{prop}\label{prop_transfaible}
   Soit $\k$ un $\QM$-anneau et $\l\in\k$. Si $X_\l(\k)\neq\vide$, l'action 
de $\exp(\xG_0)$ est transitive sur chaque $\tr{X}{n}_\l(\k)$ et sur 
$X_\l(\k)$. 
\end{prop} 
\begin{preuve}
   On raisonne par récurrence, le cas $n=0$ étant trivial.  
Supposons le résultat établi pour un entier 
$n$. Tout d'abord, 
$\tr{X}{n+1}(\k)$ contient $\tr{\pi}{n+1}(X(\k))$ et n'est donc pas vide. 
Soient $\tr{\Phi}{n+1}_1$ et $\tr{\Phi}{n+1}_2$ deux éléments de 
$\tr{X}{n+1}_\l(\k)$. Notons $\tr{\Phi}{n}_1$ et $\tr{\Phi}{n}_2$ leurs 
images par $\tr{\pi}{n}$. Par hypothèse de récurrence, il existe 
$\psi\in\xG_0(\k)$ tel que
$ \tr{\Phi}{n}_2 = \tr{\pi}{n}\exp(s_\psi)\tr{\Phi}{n}_1$.

Soit $\Psi$ l'image 
de $\tr{\Phi}{n+1}_1$ par $\tr{\pi}{n+1}\exp(s_\psi)$. C'est un élément 
de $\tr{X}{n+1}_\l(\k)$, dont l'image par $\tr{\pi}{n}$ est $\tr{\Phi}{n}_2$.
D'après la proposition \ref{prop_etape}, 
$\tr{\Phi}{n+1}_2-\Psi$ est un élément de $\xG_0(\k)$, 
homogène de degré $n+1$, donc de la forme $s_{\psi_{n+1}}(1)$. 
Comme $\tr{\pi}{n+1}\exp(s_{\psi_{n+1}})$ est précisément  
l'addition de $s_{\psi_{n+1}}(1)$, ceci achève la récurrence.
Le passage à la limite ne pose pas de problème. 
\end{preuve}

Comme $1$ est élément de $X_0(\QM)$, on a 
déjà l'isomorphisme de $\exp(\xG_0)$ et de $X_0$.

\begin{prop} \label{prop_existence}
Il existe un élément de $X_1(\QM)$.
\end{prop}
\begin{preuve}
   Il suffit de montrer pour tout $n$ que tout $\Phi\in\tr{X}{n}(\QM)$ 
peut se relever à $\tr{X}{n+1}(\QM)$. De l'hypothèse v), on 
déduit par l'action de $\Mm$ un élément $\Psi$ de $X_1(\KM)$. 
La proposition précédente fournit $\psi\in\xG_0(\KM)$ tel que 
$\tr{\pi}{n}\exp(s_\psi)$ envoie $\tr{\pi}{n}(\Psi)$ sur $\Phi$. 
L'image de $\tr{\pi}{n+1}(\Psi)$ par $\tr{\pi}{n+1}\exp(s_\psi)$
est un relevé à coefficients dans $\KM$ de $\Phi$. On conclut avec
la proposition \ref{prop_etape}.
\end{preuve}

Pour achever la preuve de la proposition \ref{prop_trans}, 
il reste à vérifier 
que $X_\l(\k)$ n'est jamais vide. Il suffit pour cela de considérer
l'action homogène de $\l$ sur un élément de $X_1(\QM)$.  

\subsection{Preuve du théorème \ref{theo_moi}}
Si $\G$ est de cardinal au moins 3, la proposition \ref{prop_trans} 
s'applique 
directement à $X=\DMRP$ (resp. $\DMRD$). On prend l'action 
de $\AM^1$ sur $\sernc{\k}{\XB}$ définie par le poids, pour laquelle 
l'hypothèse de finitude est déjà vraie. 
La stabilité de $X$ pour l'action de $\AM^1$ est laissée au lecteur. 
Les propriétés i) et ii) et l'homogénéité dans iii) 
sont évidentes. 
La propriété iii) est l'énoncé de la proposition 
\ref{prop_DMRstable} ; 
iv) est vraie par construction de $\MT$ ; v) est donnée par $\ItC$,
élément de $\DMRD_{2i\pi}(\CM)$.
 
Dans les cas $\G=\triv$ et $\G=\{\pm 1\}$ apparaît une difficulté.
La flèche $\DMRP\to\AM^1$, donnée dans ces cas par le coefficient 
de $y_{2,1}$, est de degré 2:
l'action homogène de $\mu\in\k$ 
sur $\DMRP_\l(\k)$ est à valeurs dans $\DMRP_{\l\mu^2}(\k)$. 

Tous les arguments utilisés dans \ref{subsec_approx} restent valables, 
à l'exception de la dernière phrase (pour la proposition 
\ref{prop_existence} remarquer que les racines carrées existent dans $\KM$, 
qui est ici $\CM$).
La conclusion reste vraie si l'on prouve l'existence d'un élément 
pair de $\Phi\in\DMRD_1(\QM)$, \ie{} dont toutes les 
composantes de degré impair sont nulles: Dans $\k[\sqrt{\l}]$,
l'action homogène de $\sqrt{\l}$ sur $\Phi$ ne fait intervenir que des
puisances paires de $\sqrt{\l}$ et fournit donc un élément de $\DMRD_\l(\k)$.

Soit $\Phi\in\DMRD_1(\QM)$. L'image de $\Phi$ par l'action homogène de $-1$, 
étant encore dans $\DMRD_1(\QM)$, est de la forme $\exp(s_\psi)(\Phi)$, avec
$\psi\in\dmrd_0(\QM)$. On vérifie facilement que $\exp(s_\psi/2)(\Phi)$ 
est stable par l'action homogène de $-1$, et est donc pair.

\Subsection{Conséquences}\label{subsec_conseq}
\sss{Théorème d'Écalle} Fixons un plongement $\iota : \G\to\CM^*$. 
Un élément $\Phi$ (pair si $\G=\triv$ ou $\G=\{\pm1\}$) 
de $\DMRD_1(\QM)$ fournit un isomorphisme de schémas 
$$ \AM^{1}\times\dmrd_0\longisomto\DMRD $$
qui n'est pas canonique. L'algèbre affine de $\DMRD$, 
qu'on peut voir comme engendrée par des symboles formels représentant les 
valeurs de polylogarithmes multiples sur $\G$ et soumis aux relations DMRD, 
est donc isomorphe au produit tensoriel de $\QM[t]$ et de l'algèbre
symétrique formée sur le dual gradué de $\dmrd$.\footnote{Cette dernière est
aussi la duale graduée de l'algèbre enveloppante universelle de $\dmrd$.} L'isomorphisme respecte
les graduations, à condition d'attribuer le degré 2 à $t$ si $\G=\triv$ 
ou $\G=\{\pm 1\}$. Un résultat analogue vaut pour $\DMRP$. 

Pour $\G=\triv$, ceci est le théorème d'Écalle \cite{Ecalle}, dont la 
démonstration repose en partie sur des méthodes similaires (linéarisation
par une algèbre de Lie isomorphe à $\dmrd$, \cf{} 
\cite[Appendice A]{thesejoe}).   
\sss{Irréductibles de Drinfel'd} 
Les propriétés énumérées en \ref{sss_proj} s'appliquent en particulier 
aux automorphismes de $\G$, qui fournissent donc par image directe des 
automorphismes de $\DMRD, \dmrd$, etc. On note dans la suite 
$N$ l'ordre de $\G$ et $\phi$ la fonction indicatrice d'Euler.

La composition $F$ de l'action homogène de $-1$ avec l'application 
image directe associée à $\s\mapsto\s^{-1}$ est 
la substitution $x_\s\mapsto -x_{\s^{-1}}, x_0\mapsto-x_0$. 
Elle est involutive. C'est, pour notre 
cas, ce que Deligne appelle le Frobenius réel \cite{DelPi1}.

Comme $F$ fixe les $x_\s-x_{\s^{-1}}$, elle commute à $\DMRP\to\AM^1$. 
Par le théorème \ref{theo_moi}, 
il existe donc un élément $\psi$ de $\dmrd(\CM)$, tel que 
$F(\ItC)=\exp(-s_\psi)(\ItC)$. De $F^2=\id$, on déduit $F(\psi)=-\psi$. 
Cet élément dépend du plongement $\iota$ 
de $\G$ dans $\CM^*$. On obtient ainsi 
$\phi(N)$ éléments $\psi_\iota$. Ils sont 
permutés par les automorphismes de $\G$ et égaux au signe près 
à leurs images par l'inversion. On les indice par le choix d'un 
$\iota$ dans chaque orbite de l'inversion. 

Soit $\psi_{n,\iota}$ la composante homogène de degré $n$ de $\psi_\iota$. 
Dans le cas $N=1$, les $\psi_n$ sont les éléments irréductibles de 
$\grt_1(\CM)$ exhibés par Drinfel'd \cite[p. 860]{DrinQTQH}.
D'après \cite{Del01}, les composantes homogènes non nulles des 
$\psi_{n,\iota}$ 
engendrent en général l'image dans $\mt$ de l'algèbre de Lie du  
groupe motivique $\UC_\gr$. Appliquant la propriété (M) de l'introduction 
dans les deux cas,
on obtient que la variété des relations d'origine motivique est incluse 
dans $\DMRD$. Autrement dit, les relations DMRD sont d'origine motivique. 

On calcule facilement le terme de longueur 1 des $\psi_{n,\iota}$, 
 en utilisant \ref{ssslng1}: 

Dans le cas réel $N\leq 2$, l'inversion de $\G$ est l'identité 
et $F$ se réduit à l'action homogène de $-1$. La série $\psi$ est 
impaire. Pour $n\geq 1$, son terme de longueur 1 et de poids $n$ 
vaut $2\z(2n+1)y_{2n+1,1}$. 
On a $\psi_1=0$ pour $N=1$ et $\psi_1=2\log(2)y_{1,-1}$ pour 
$N=-1$. 

Dans le cas général, $\psi_{n,\iota}$ vaut $\sum\limits_{\s\in\G} (L_n(\iota(\s)) - (-1)^n L_n(\iota(\s^{-1})))y_{n,\s} $

Les $\psi_{n,\iota}$ non nuls ont donc 
un terme de longueur non nul. Ils sont linéairement indépendants, et 
irréductibles, car le crochet d'Ihara est homogène pour la longueur.

\begin{prob}
Pour quelles valeurs de $N$ les $\psi_{n,\iota}$ 
engendrent-ils $\dmrd_0$ ? Sont-ils libres ? 
\end{prob}
Comme mentionné dans l'introduction, on n'attend pas en général 
de réponse positive. Pour $N=1$, ce problème est une variante des 
questions de Drinfel'd \cite{DrinQTQH} à propos de $\grt_1$, parfois 
qualifiée de \og{}conjecture de Deligne-Drinfel'd.\fg{} La première 
question de la variante pour 
l'image de l'algèbre de Lie du complété pro-$\ell$ de
$\Gal(\ov{\QM}/\QM)$ dans $\grt_1(\QM_\ell)$ (conjecture de Deligne-Ihara) 
a été résolue par Hain et Matsumoto \cite{HainMat}. 

L'algèbre de Lie de $\UC_\gr$ est libre, engendrée par des éléments 
dont les images dans $\mt$ sont les $\psi_n$. La première question 
revient donc à demander si toutes les relations d'origine motivique
proviennent de DMRD. La seconde est équivalente à l'injectivité de 
$\UC_\gr\to\MT$.
 
Deligne a prouvé la liberté pour $N\in\{2,3,4\}$ (non publié). 
Pour $N=1$, on a obtenu par ordinateur une réponse positive
aux deux questions jusqu'en poids 19; 
cela fera l'objet d'un autre article \cite{ENR}. 

Les conjectures de transcendance et la conjecture de Deligne-Drinfel'd amènent 
également à la question transversale:
\begin{prob}
   Les algèbres de Lie $\grt_1$ et $\dmr_0$ sont elles égales ? 	
\end{prob}
Une réponse positive ramènerait la construction 
explicite d'associateurs rationnels à celle d'éléments de $\DMRP_1(\QM)$,
ce qu'on espère être plus facile.

\nocite{IhICM, Hoffalg, Demaz, IhIsrael, Broad, Chen}

\bibliographystyle{joeplain}
\bibliography{defs,motifs,quantique,divers,homotopie,polychoses,combi,moi,traites,ICM}
\setcounter{tocdepth}{2}
\newpage
\tableofcontents
\end{document}